\documentclass[12pt,a4paper]{article}
\usepackage{amssymb}
\usepackage{amsfonts}
\usepackage{amsmath}
\usepackage{amsthm}
\usepackage[numbers]{natbib}
\usepackage{setspace}

\usepackage{caption}
\theoremstyle{definition}

\newtheorem{rem}{Remark}
\newtheorem{ex}{Example}
\theoremstyle{plain}
\newtheorem{thm}{Theorem} 
\newtheorem{lem}{Lemma}
\newtheorem{cor}{Corollary}
\newtheorem{prop}{Proposition}
\newtheorem{claim}{Claim}
\newenvironment{pf}[1][Proof]{\noindent\textbf{#1.} }{\ \rule{0.6em}{0.6em}}
  
\usepackage[normalem]{ulem}
\usepackage{color}
\definecolor{shade}{rgb}{0.9,0.3,0.3}

\newcommand{\eps}{\varepsilon}
\newcommand{\bs}{\backslash}

\newcommand{\R}{\mathbb{R}}
\newcommand{\N}{\mathbb{N}}

\begin{document}

\title{No-regret Dynamics and Fictitious Play
\thanks{The authors thank William Sandholm whose comments led to a substantial
improvement of the paper, as well as Mathieu Faure, Sergiu Hart, Josef Hofbauer, Alexander Matros, Karl Schlag, Eilon Solan and Sylvain Sorin for helpful comments and suggestions.}}

\author{Yannick Viossat\thanks{CEREMADE, Universit{\'e} Paris-Dauphine, Place du Mar{\'e}chal de Lattre de Tassigny, F-75775 Paris, France. \emph{E-mail:} viossat $\alpha\tau$ ceremade.dauphine.fr}
\and Andriy Zapechelnyuk
\thanks{School of Economics and Finance, Queen Mary, University of London,
Mile End Road, London E1 4NS, UK. \emph{E-mail:} a.zapechelnyuk $\alpha\tau$ qmul.ac.uk
} 
}
\date{September 5, 2012}

\maketitle
\begin{abstract}
\begin{spacing}{1.2}
\noindent Potential based no-regret dynamics are shown to be related to fictitious play. Roughly, these are $\eps$-best reply dynamics where $\eps$ is the maximal regret, which vanishes with time. This allows for alternative and sometimes much shorter proofs of known results on convergence of no-regret dynamics to the set of Nash equilibria. 


\end{spacing}
\end{abstract}

\newpage
\section{Introduction}

\emph{No-regret strategies} are simple adaptive learning rules that
recently received a lot of attention in the literature.\footnote{These rules have been used to investigate convergence to equilibria in the context of learning in games \citep{Freund99, Fudenberg95, Hart00, Hart01Gen, Hart03}, for combining different forecasts \citep{Foster93,Foster99} (for an overview of the forecast combination literature see \cite{Clemen07,Timmerman06}) and for combining opinions, which is also of interest to management science \citep{Larrick06}. In finance this method has been used to derive bounds on the prices of financial instruments \citep{Chen10,DeMarzo06}. This method can be applied to various tasks in computer science, such as job scheduling \citep{Mansour10} and routing \citep{Blum06} (for a survey of applicable problems in computer science see \citep{Irani96}).}
In a repeated game, a player has a \emph{regret} for an action if,
loosely speaking, she could have obtained a greater average payoff
had she played that action more often in the past. In the course of
the game, the player reinforces actions that she regrets not having
played enough, for instance, by choosing next action with
probability proportional to the regret for that action, as in Hart
and Mas-Colell's \cite{Hart00} \emph{regret matching} rule. Existence of {\it no-regret strategies} (i.e., strategies that guarantee no regrets almost surely in the long run) is known since \citet{Hannan57};
wide classes of no-regret strategies are identified by \citet{Hart01Gen}
and \citet{CesaBianchi03}.%
\footnote{This paper deals with the simplest notion of regret known as \emph{unconditional} (or \emph{external})
regret \citep{Fudenberg95,Hart01Gen,Hart03}. For more sophisticated
regret notions, see \citet{Hart00}, \citet{Lehrer03}, and \citet{CesaBianchi06}.}

A \emph{no-regret dynamics} is a stochastic process that
describes trajectories of the average correlated play of players and
that emerges when every player follows a no-regret strategy (different
players may play different strategies). By definition, it converges to the Hannan set (the set of all correlated actions that satisfy the no-regret condition first stated
by \citet{Hannan57}).\footnote{\label{Foot:Hannan}The Hannan set of a game is also known as the set of \emph{weak correlated equilibria} \citep{Moulin78} or \emph{coarse
correlated equilibria} \citep[Ch.3]{Young04}.} This set is typically large. It contains the set of correlated equilibria of the game and we show that it may even contain correlated actions that put positive weight \emph{only} on strictly dominated actions. Thus convergence of the average play to the Hannan set often provides very little information about what the players will actually play, as it does not even imply exclusion of strictly dominated actions.

In this paper we show that no-regret dynamics are intimately 
linked to the classical fictitious play process \citep{Brown51}. Drawing on \citet{Monderer97}, we first show that contrary to the standard, discrete-time version, continuous fictitious play leads to no regret. We then show that, for a large class of no-regret dynamics, if a player's maximal regret is $\eps >0$, then she plays an $\eps$-best reply to the average correlated play of the others. Since in this class the maximal regret vanishes (see Corollary \ref{Cor:1} below), it follows that, for a good choice of behavior when all regrets are negative, the dynamics is a vanishingly perturbed version of fictitious play. 

For two-player finite games, this observation and the theory of perturbed differential inclusions \citep{Benaim05, Benaim06} 
allow us to 
relate formally the asymptotic behavior of no-regret dynamics and of continuous fictitious play (or its time-rescaled version, the 
best-reply dynamics \citep{Gilboa91}). In classes of games in which the behavior of continuous fictitious play is well known, this provides 
substantial information on the asymptotic behavior of no-regret dynamics. In particular, we recover most known convergence 
properties of no-regret dynamics. Our results do not just 
allow us to find new and sometimes much shorter proofs of convergence of no-regret dynamics towards the set of Nash equilibria 
in some classes of games, such as dominance solvable game or potential games. They also allow us to relate the asymptotic behavior 
of no-regret dynamics and continuous fictitious play in case of divergence, as in the famous Shapley game \citep{Shapley64}.

These results extend only partially to $n$-player games (though they fully extend to $n$-player games with linear incentives \citep{Selten95}). The issue is that in $n$-player games no-regret dynamics turn out to be related to the correlated version of continuous fictitious play, in which the players play a best-reply to the \emph{correlated} past play of the others. This version of fictitious play is defined through a correspondence which is not convex valued. This creates technical difficulties, because the theory of perturbed differential inclusions is not developed for non convex valued correspondences.

A different way to analyze no-regret dynamics is to show that some sets attract nearby solution trajectories. We show that strict Nash equilibria and, more
generally, the intersection of the Hannan set and the sets that are \emph{closed under rational behavior (curb)}\footnote{A product set of action profiles is called \emph{closed under rational
behavior (curb)} \citep{Basu91} if it contains all best replies of each player
whenever she believes that no actions outside this
set are being played by the other players.} 
are attracting for no-regret dynamics, in a sense to be defined in Section \ref {Sec:curb}.

The remainder of the note is organized as follows. The next section introduces no-regret dynamics. Section \ref{sec:main} studies the links between no-regret dynamics and fictitious play. Section \ref{Sec:curb} shows that the intersection of the Hannan set and curb sets is attracting for no-regret dynamics. Section \ref{sec:cont} studies the continuous-time version and the expected version of no-regret dynamics. Finally, the Appendix contains the proofs of the main results, as well as counterexamples illustrating the complexity of the relationship between ICT and limit sets. 

\section{Preliminaries}

\label{Sec:Prelim}
Consider a bimatrix game $\Gamma=(A_{i},u_{i})_{i=1,2}$, where $A_{i}$
is the set of actions of player $i$ and $u_{i}:A\rightarrow\R$ 
is her payoff function, with $A=A_{1}\times A_{2}$. For any finite set $B$, denote by $\Delta(B)$ the set of probability distributions over $B$. 
A mixed action of player $i$ is an element of $\Delta(A_i)$. A correlated action $z$ is a probability distribution over the set of pure action profiles, i.e., $z \in \Delta(A)$. Given such a $z$, let $z_i \in \Delta(A_i)$ and $z_{-i} \in \Delta(A_{-i})$ denote its marginals for player $i$ and her opponent, respectively. Thus, $z_i(a_i)=\sum_{a_{-i} \in A_{-i}} z(a_i,a_{-i})$. Throughout, $-i$ refers to $i$'s opponent. As usual, let $u_{i}(z)=\sum_{a\in A}z(a)u_{i}(a)$ and $u_i(k,z_{-i})=\sum_{a_{-i}\in A_{-i}} z_{-i}(a_{-i}) u_{i}(k,a_{-i})$ for $k \in A_i$. Depending on the context, $a_i$ may refer to a pure action -- an element of $A_i$ -- or to a vertex of $\Delta(A_i)$, i.e., a Dirac measure on a pure action.

The game is played repeatedly in discrete time periods $t\in\N^*=\{1,2,\ldots\}$.
In every period $t$ each player $i$ chooses an action $a_{i}(t)\in A_{i}$
and receives payoff $u_{i}(a(t))$ where $a(t)=(a_{1}(t),a_{2}(t))$.
Denote by $h(t)=(a(1),a(2),\ldots,a(t))$ the history of play up to
$t$, and let $\mathcal{H}$ be the set of all finite histories (including
the empty history). A strategy of player $i$ is a function $q_{i}:\mathcal{H}\rightarrow\Delta(A_{i})$ that stipulates to play
in every period $t=1,2,\ldots$ a mixed action $q_{i}(t)\equiv q_{i}(h(t-1))$
as a function of the history before $t$. The weight that this mixed action puts on action $k \in A_i$ is denoted by $q_{i,k}(t)$ 

The \emph{average correlated play} up to period $t$ is $z(t) = \frac{1}{t} \sum_{\tau=1}^t a(\tau)$, where we identify $a(\tau)$ with the corresponding vertex of $\Delta(A)$. Since $z(t)=\frac{1}{t} \left[ a(t)+(t-1)z(t-1)\right]$, it follows that for all $t >1$, 
\begin{equation}
z(t)-z(t-1)=\frac{1}{t}\left(a(t)-z(t-1)\right). \label{Eq:DT}
\end{equation}

For a correlated action $z$, the \emph{regret} of player $i$ for action $k$ is defined as $R_{i,k}(z)= u_i(k,z_{-i}) - u_i(z)$, and her maximal regret as $R_{i,\max}(z)=\max_{k \in A_i} R_{i,k}(z)$. Typically we deal with the regret based on the average correlated play, $z(t)$, up to some period $t$. In this case the regret of player $i$ for action $k \in A_i$ is equal to the difference between the average payoff she would have obtained by always playing $k$ (assuming that her opponent's play remains the same) and her average realized payoff:
\[
R_{i,k}(z(t)) = u_i(k,z_{-i}(t)) - u_i(z(t))= \frac{1}{t} \sum_{\tau=1}^t [u_i(k,a_{-i}(\tau)) - u_i(a(\tau))].
\]
To simplify notations, we will often write $R_{i,k}(t)$ for $R_{i,k}(z(t))$ and $R_{i,\max}(t)$ for $R_{i,\max}(z(t))$. 

Player $i$ has no asymptotic regret if her average realized payoff is asymptotically no less than her best-reply payoff against the empirical distribution of her opponent:  
\begin{equation}
\limsup_{t\rightarrow\infty} R_{i,\max}(t) \leq 0.\label{Eq:Obj}
\end{equation}
A strategy of player $i$ is a \emph{no-regret strategy} if for any strategy of the other player, inequality (\ref{Eq:Obj}) holds almost surely.  This property is also called \emph{Hannan consistency} \citep{Hart01Gen} or \emph{universal consistency} \citep{Fudenberg95}.

It is well known in the literature since \citet{Hannan57}
that there exist simple no-regret strategies. \citet{Hart01Gen} describe a wide class of \emph{potential
based} no-regret strategies. A twice differentiable, convex function $P_{i}:\R^{A_{i}}\rightarrow\R$
is called a \emph{potential} if it satisfies the following conditions: 
\begin{description}
\item [{(R1)}] $P_{i}(\cdot)\geq0$, and $P_{i}(x)=0$ for all $x\in\R_{-}^{A_{i}}$; 
\item [{(R2)}] $\nabla P_i(\cdot) \geq 0$, and $\nabla P_i(x) \cdot x >0$ for all $x \notin \R^{A_i}_{-}$;
\item [{(R3)}] if $x \notin \R^{A_i}_{-}$ and $x_{k}\le 0$, then $\nabla_k P_{i}(x)=0$,
\end{description}
where $\nabla_k$ denotes the partial derivative with respect to $x_i(k)$. 
The potential $P_{i}$ can be viewed as a generalized distance
function between a vector $x\in\R^{A_{i}}$ and the nonpositive orthant $\R_{-}^{A_{i}}$. Let $R_i(t)=(R_{i,k}(t))_{k \in A_i}$ denote player $i$'s regret vector. 
\begin{prop}
\label{Thm:HM} Let $P_{i}$ satisfy (R1)--(R3) and let 
strategy $q_{i}$ satisfy 
\begin{equation*}\label{Eq:PB}
q_{i,k}(t+1)=\frac{\nabla_k P_{i}(R_{i}(t))}{\sum_{s\in A_{i}}\nabla_{s} P_{i}(R_{i}(t))},\text{ \ } \forall k\in A_{i}, \tag{Q1}
\end{equation*}
whenever $R_{i,\max} (t) > 0$. Then $q_{i}$
is a no-regret strategy.
\end{prop}
\begin{pf}
This holds by Theorem 3.3 of \citet{Hart01Gen}, whose conditions (R1) and (R2) are satisfied by our conditions (R1)--\eqref{Eq:PB} and (R2), respectively, and whose proof is based on the Blackwell's Approachability Theorem \citep{Blackwell56}.
\end{pf}

A standard example of no-regret strategy satisfying the above conditions is obtained by letting $P_i$ be the $l_{\mathbf{p}}$-norm
on $\R_{+}^{A_{i}}$, i.e. $P_i(x)=(\sum_{k\in A_{i}}[x_{k}]_{+}^{\mathbf{p}})^{1/\mathbf{p}}$
with $1<\mathbf{p}<\infty$, where $[x_{k}]_{+}=\max(0, x_k)$. The resulting strategy $q_{i}$ is called the $l_{\mathbf{p}}$\emph{-norm}
strategy \citep{CesaBianchi03, Hart01Gen}. It is defined by
\[
q_{i,k}(t+1)=\frac{[R_{i,k}(t)]_{+}^{\mathbf{p}-1}}{\sum_{s\in A_{i}}[R_{i,s}(t)]_{+}^{\mathbf{p}-1}}, \text{ \ } \forall k\in A_{i},\]
 whenever $R_{i,\max} (t) > 0$. The $l_{2}$-norm strategy is the \emph{regret-matching strategy} \citep{Hart00}, that stipulates to
play an action in the next period with probability proportional to
the regret for that action. For large $\mathbf{p}$, the $l_{\mathbf{p}}$-norm strategies approximate
fictitious play. 

We say that the average correlated play $z(t)$ follows a \emph{no-regret dynamics} if both players use (possibly different) no-regret strategies. A trajectory $(z(t))_{1 \leq t \leq +\infty}$ of a no-regret dynamics is thus a solution of (\ref{Eq:DT}) where $a(t)$ is a realization of $(q_1(t), q_2(t))$ and $q_1, q_2$ are no-regret strategies.  We focus on the class $\mathcal{R}$ of no-regret dynamics such that: 

(i)  the no-regret strategies $q_1$, $q_2$ of the players are potential-based: they satisfy (Q1) for some potentials $P_1$, $P_2$ satisfying (R1)-(R3);

(ii) if a player has no-regret  then he takes some constant pure action: for each $i=1,2$, there exists $c \in A_i$ such that 
\begin{equation*}\label{Eq:PBR}
a_i(t+1) = c \ \ \text{whenever $R_{i,\max}(t)\le 0$}. \tag{Q2}
\end{equation*}
Our results are valid for a somewhat wider class of no-regret dynamics. What we really need, beside a no-regret dynamics, is that from some period $t_0$ on: 

(i$'$) if a player has positive regret for some actions, then she plays one of these actions. 

(ii$'$) if a player never has any positive regret, then she plays an $\eps(t)$-best reply to the empirical distribution of her opponent, 
where $\eps(t)=\eps(h(t)) \to 0$ almost surely.

\begin{rem}
\label{rem:Q1}
Property (i$'$) follows from (R3) and (\ref{Eq:PB}). This is a \emph{better reply property} that stipulates to assign a positive probability only on better reply actions to the opponent's empirical distribution of play (``better" with respect to the realized payoff). Also it implies that if $R_{i,\max}(t) >0$ in some period $t$, then $R_{i,\max}(t')>0$ for all $t' > t$. 
Indeed, when an action $k$ with positive regret is played, the 
sign of $R_{i,k}(t)$  does not change, hence the maximal regret remains positive \citep[Proposition 4.3]{Hart01Gen}. 
\end{rem}     
\begin{rem} \label{rem:Q2} 
Assumption (Q2) is a simple way of ensuring (ii$'$), and in addition, that if $R_{i,\max}(t) \leq 0$ for all $t$, then $R_{i,\max}(t) \to 0$ as $t \to +\infty$.\footnote{This additional property is needed for Corollary \ref{Cor:1} below, but for our main results (ii$'$) suffices.}  Indeed, if $R_{i,\max}(t) \leq 0$ for all $t > t_0$ then by (Q2), 
for all $t > t_0$, $t R_{i,c}(t)=t_0 R_{i,c}(t_0)$, hence $R_{i,c}(t)\to 0$. It follows that 
$R_{i,\max}(t) \to 0$ and that for all $t>t_0$, player $i$ plays an $\eps(t)$-best reply with 
$\eps(t):= \max_{k \in A_i} u_{i}(k, z_{-i}(t)) - u_i(c, z_{-i}(t))= R_{i,\max}(t) - R_{i,c}(t) \to 0$. For a discussion of other possible assumptions, see \citet{Hart03}, Appendix A.
\end{rem}

Note that there are no-regret dynamics that do not satisfy (i$'$). For instance, stochastic fictitious play with a noise parameter that declines with time at an appropriate rate (see, e.g., \citet{Benaim11}). This process is not potential based in our sense due to the time inhomogeneity, but this is not the crucial point, since (i$'$)-(ii$'$) would suffice.

Define the \emph{Hannan set} $H$ of the stage game $\Gamma$ as the set of all correlated actions of the players where each player has no regret: 
\[
H=\left\{ z\in\Delta(A)\,\left|\,\max_{k\in A_{i}}u_{i}(k,z_{-i})\leq u_{i}(z)\,\mbox{ \ for each }i=1,2\right.\right\} .
\]
The {\it reduced Hannan set} $H_R$ is the subset of  $H$ in which at least one regret is exactly zero for each player:
\[
H_R=\left\{ z\in\Delta(A)\,\left|\,\max_{k\in A_{i}}u_{i}(k,z_{-i})= u_{i}(z)\,\mbox{ \ for each }i=1,2\right.\right\} .
\]

The next property of no-regret dynamics is straightforward by the definition of no-regret strategies and Remark \ref{rem:Q2} (see, e.g., \citet[Corollary 3.2]{Hart03}). 
\begin{cor}
\label{Cor:1} For every no-regret dynamics in class $\mathcal{R}$, the trajectories converge almost surely to the reduced Hannan set.
\end{cor}
Convergence of the average play $z(t)$ to set $H_R$ does not
imply its convergence to any particular point in $H_R$. Moreover, even if $z(t)$ converges
to a point, this point need not be a Nash equilibrium.

\section{Fictitious play and no-regret dynamics} \label{sec:main}

\subsection{Fictitious play} \label{Sec:DFP}

In {\it discrete fictitious play}, in every period $t$ after the initial one, player $i$ plays a pure best reply $a_i(t)$ to the average past play of her opponent $x_{-i}(t-1):=\frac{1}{t-1} \sum_{\tau=1}^{t-1} a_{-i}(\tau)$ (here $a_{-i}(\tau)$ is a vertex of $\Delta(A_{-i})$). The latter is called the \emph{belief} of player $i$ on her opponent's next move. Formally, for any $x=(x_1,x_2)$ in $\Delta(A_1) \times \Delta(A_2)$, denote by $BR_i (x_{-i})$ player $i$'s set of best replies to $x_{-i}$:
\begin{equation*}
BR_i (x_{-i}):= \Big\{x_i \in \Delta(A_i) \, \Big| \, u_i(x_i, x_{- i})  = \max_{k\in A_i} u_i(k,x_{-i})\Big\}, \ \ i=1,2.
\end{equation*}
Let $BR(x)=BR_1(x_{2})\times BR_2(x_{1})$. A discrete-time trajectory $(x(t))_{t=1}^{\infty}$ on $\Delta(A_1)\times\Delta(A_2)$ is a solution of {\it discrete fictitious play} (DFP) if for every $t>1$

\begin{equation}\label{Eq:FP}
x(t)-x(t-1)= \frac{1}{t}\left(a(t)-x(t-1)\right)
\end{equation}
where $a(t)=(a_1(t),a_2(t))$ and $a_i(t) \in BR_i(x_{-i}(t-1))$ is a vertex of $\Delta(A_i)$ associated with some pure best reply action, $i=1,2$.   

Analogously, an absolutely continuous function $x:[1,\infty)\to \Delta(A_1)\times\Delta(A_2)$ is a solution of {\it continuous fictitious play} (CFP) if for almost all $t \geq 1$, $x(t)$ is differentiable and 
\begin{equation*}
\dot{x}(t)= \frac{1}{t}\left(q(t)-x(t)\right),
\end{equation*}
where $q(t) \in BR(x(t))$ is now a profile of \emph{mixed} actions. This may be written as the differential inclusion:  
\begin{equation}\label{Eq:CFP2}
\dot{x}(t) \in \frac{1}{t}\left(BR (x(t)) -x(t) \right).
\end{equation}
The average correlated play satisfies $z(t):=\frac{1}{t} \left(z(1) + \int_1^{t} q(\tau) d\tau \right)$ for some initial condition $z(1)$ such that $z_i(1)=x_i(1)$, $i=1,2$. Thus, for almost all $t$, $z(t)$ is differentiable and   
\begin{equation} \label{eq:cont}
\dot z(t) = \frac{1}{t}(\bar{q}(t)- z(t)),
\end{equation}  
where $\bar{q}=q_1 \otimes q_2 \in \Delta(A)$ is the product distribution corresponding to the mixed strategy profile $q=(q_1,q_2) \in \Delta(A_1) \times \Delta(A_2)$, and $q_i$ is a best-reply to $z_{-i}$.\footnote{\label{ft:BRD}This definition of CFP guarantees that solutions exist in all games and for all initial conditions, and that by the change of time scale $y(t)=x(e^t)$, CFP corresponds to the {\it best-reply dynamics} \citep{Gilboa91, Matsui92} defined by $\dot{y} \in BR (y) - y$.
Another definition of CFP (e.g., \citet[p.~445]{Monderer97} and \citet[pp.~252--253]{Berger07}) consider only trajectories that are piecewise linear, such that $q_i(t)$ is always a pure action (technically, a vertex of $\Delta(A_i)$), and that the times at which $q(t)$ changes have no finite accumulation point. This restricted definition is easier to handle, but in many games there do not exist such trajectories from every initial condition.
}

In discrete or continuous fictitious play, the marginals $z_1(t)$, $z_2(t)$ of the average past play are equal to the beliefs $x_1(t)$, $x_2(t)$. By analogy, if $z(t)$ is the average past play generated by a no-regret dynamics, it is convenient to call $z_{-i}(t)$ the belief of player $i$ about her opponent's next move. This illuminates a crucial difference between fictitious play and no-regret dynamics in class $\mathcal{R}$: under fictitious play, a player chooses a {\it best reply} to her belief, whereas under no-regret dynamics, she chooses a {\it better reply} (``better" with respect to her average realized payoff). 

\subsection{Continuous fictitious play leads to no regret}

It is well known that discrete fictitious play does not lead to {\it no regret} \citep{Hart01Gen,Young93E}. Consider the following example:
\begin{figure}[!htb]
\[
\begin{array}{c|cc}
& L & R \\
\hline
L & 1, \sqrt{2}  & 0,0 \\
R & 0, 0 & \sqrt{2}, 1
\end{array}
\]
\caption{}\label{F:DFP}
\end{figure}

\noindent Because $\sqrt{2}$ is irrational, $L$ and $R$ cannot both be best-replies to the empirical past play of the other player. Thus, any DFP process is entirely determined by its first move. Assume that the first move is off the diagonal, say $(L,R)$. Due to the symmetry of the game and the absence of ties, both players always switch to another action simultaneously. Therefore the play is locked off the diagonal and the maximal regret is at least $\sqrt{2}/ (1+\sqrt{2})$ at any stage. This holds in the mixed extension of the game, since at any stage the players have a unique, pure best reply.

Since the continuous fictitious play process is a continuous-time version of DFP, intuitively, it should not lead to {\it no regret} either. The following result --- a generalization of Theorem D of \citet{Monderer97} --- shows that this intuition is misleading.

\begin{prop}\label{P:CFP}
Under any solution of continuous fictitious play, the average correlated play converges to the reduced Hannan set.
\end{prop}

This discrepancy between DFP and CFP may be explained as follows. Playing an action with positive regret decreases the regret for this action. In CFP, roughly, when an action is played it remains a best reply, hence it is associated with maximal regret for a small time increment. Precisely, the derivative of the regret for the action played is equal to the derivative of the maximal regret. Since the regret for this action decreases, so does the maximal regret. In contrast, in DFP, an action played at stage $t$ has maximal regret at stage $t$, but not necessarily at stage $t+1$. Thus the fact that the regret for this action decreases does not entail that the maximal regret does.

\begin{pf}[Proof of Proposition \ref{P:CFP}]
For comparison with \citet[Theorem 3.1]{Hart03}, rescale time (let $\tilde{t}=\exp t$) so that (\ref{eq:cont}) becomes $\dot z = \bar{q}-z$. For any mixed action $\sigma_i \in \Delta(A_i)$, let
\[R_{i,{\sigma_i}(t)}:=\sum_{k \in A_i} \sigma_i(k) R_{i,k}(t)= u_i(\sigma_i, z_{-i}(t)) - u_i(z(t))\]   
Let $v_i(t)=R_{i,\max}(t)$. Note that $R_{i,k}$ is Lipschitz continuous for all $k$ in $A_i$. Thus it follows from Theorem A.4 of \citet{Hofbauer09JET} that, 
for almost all $t$, $v_i$ and $R_{i,k}$ are differentiable, and for all $k$ such that  $q_{i,k}(t)>0$, we have $\dot{v}_i(t)=\dot{R}_{i,k}(t)$. It follows that $\dot{v}_i=\sum_k q_{i,k} \dot{R}_{i,k}= \dot{R}_{i,q_i}$. Furthermore:
\[ \dot{R}_{i,q_i} = u_{i}(q_i, \dot{z}_{-i})-u_{i}(\dot{z})=u_{i}(q_i, q_{-i}-z_{-i})-u_{i}(\bar q-z)
= -[u_{i}(q_i,z_{-i})-u_{i}(z)] = -R_{i,q_i}= - v_i.\label{Eq:A-1b}\]
Thus, $\dot{v_i}=-v_i$. Therefore, $v_i(t)$ converges to zero for all $i=1,2$, hence $z(t) \to H_r$.
\end{pf}

\begin{rem} \label{rem:uni}
In the proof, we did not use that $q_{-i}$ is a best-reply to $z_i$.  This shows that the fact that CFP leads to no-regret is a unilateral property. That is, if a player's behavior evolves according to CFP, then she has no asymptotic regret, independently of her opponent's behavior (see also \citet[p.~445]{Monderer97}).
\end{rem} 

\begin{rem} 
CFP and the best-reply dynamics converge to the set of Nash equilibria in finite zero-sum games \citep{Hofbauer06}.  The usual proof is to show that the ``duality gap"  $W(x)=\max_{k \in A_1} u_1(k, x_2) - \min_{s \in A_2} u_1(x_1, s)$ converges to zero. This follows from the above proof, since in a two-player zero sum game $W(x(t))=R_{1,\max} (z(t)) + R_{2,\max}(z(t))$, where $x$ is a solution of CFP and $z$ the associated correlated play.
\end{rem}

\subsection{No-regret dynamics is perturbed CFP}
 
In the previous subsection we showed that CFP leads to no regret. Conversely, we now show that any no-regret dynamics in class $\mathcal R$ (as defined in Section \ref{Sec:Prelim}) is closely related to CFP. We first explain the intuition. Denote by $BR_i^{\eps} (x_{-i})$ the set of $\eps$-best replies of player $i$ to the mixed action $x_{-i}$ of her opponent: 
\begin{equation*}
BR_i^{\eps} (x_{-i})= 
\Big\{x_i \in \Delta(A_i) \ \Big| \ u_i(x_i, x_{- i})  \geq \max_{k\in A_i} u_i(k,x_{-i})-\eps \Big\}, \ \ i=1,2.
\end{equation*}
The crucial observation is the following. 
\begin{lem} \label{lm0} 
Assume that the maximal regret is less than $\eps$. Then any action with positive regret is an $\eps$-best reply to the average play of the opponent.
\end{lem}
\begin{pf}
If player $i$ has  positive regret for action $a_i$ at some $z \in \Delta(A)$, then $u_i(z) - u_i(a_i, z_{-i}) < 0$. But by assumption $\max_{k \in A_i} u_i(k, z_{-i}) - u_i(z) \leq \eps$. Therefore, $\max_{k \in A_i}$ $u_i (k, z_{-i}) - u_i(a_i, z_{-i}) < \eps$, and $a_i$  is an $\eps$-best reply to $z_{-i}$.
\end{pf}
\bigskip

Since no-regret dynamics in class $\mathcal{R}$ only pick actions with positive regret, they only pick $\eps$-best replies to the average play of the others, where $\eps$ is the maximal regret. Since this maximal regret approaches zero almost surely, eventually only almost-exact best replies are picked.
This provides the intuition why no-regret dynamics and fictitious play may exhibit similar asymptotic behavior. Finding a precise link, however, is not obvious. For instance, there could exist actions that are $\eps_t$-best replies in each period $t$, with $\eps_t \to 0$, but never exact best replies. Thus a limit play of no-regret dynamics may include such actions, but this cannot happen under fictitious play. 
\begin{figure}[!htb]
\[
\begin{array}{c|cc}
& L & R \\
\hline
L & 1, 0  & 0,\sqrt{2} \\
R & 0, 1 & \sqrt{2}, 0 \\
C & \eta, 0 & \eta, 0
\end{array}
\]
\caption{
}
\label{F:EMP}
\end{figure}

Consider the example shown on Fig.~\ref{F:EMP}.
Let $\eta=\sqrt{2}/(1+\sqrt{2})$. It is easy to verify that action $C$ is player 1's best reply to player 2's mixed action $x_2$ if and only if $x_2=(\eta,1-\eta)$. Let us first consider DFP. Since $\eta$ is an irrational number, after every finite history of play, $C\not\in BR_1(x_2(t))$; consequently DFP never picks $C$ (except, possibly, at the initial period).\footnote{Starting with an arbitrary belief $x_2(1)$ would not help since $C$ is  a best-reply only when $x_2(t)=(\eta, 1-\eta)$, which happens at most once.}  However, it may be shown that under any DFP trajectory, the average play $x_2(t)$ of player $2$ converges to $(\eta,1-\eta)$, to which $C$ is a best-reply. It follows that $C$ is an $\eps_t$-best reply to $x_2(t)$ for some sequence $\eps_t\to 0$. Thus a no-regret dynamics with the same trajectory of the marginal play of player 2 {\it might} choose action $C$ a positive fraction of time in the long run. 

This example does not apply to CFP, as in this case $x_2(t)$ need not be a rational number; and as we show below, the asymptotic behavior of no-regret dynamics and CFP can be formally related using the theory of perturbed differential inclusions \citep{Benaim05, Benaim06}.

Before stating a precise result, we need some definitions. A set $L \subset \Delta(A_1) \times \Delta(A_2)$ is \emph{invariant} under CFP if for every initial point $x \in L$ there exists a solution $x(\cdot)$ of CFP, defined for all $t  >0$ (not only $t \geq 1$) and such that $x(1)=x$ and $x(t) \in L$ for all $t >0$. A nonempty compact invariant set is an \emph{attractor} if it attracts uniformly all trajectories starting in its neighborhood. An invariant set $L$ is \emph{attractor-free} if no proper subset of $L$ is an attractor for the dynamics restricted to $L$. A nonempty compact set $L$ is \emph{internally chain transitive} (ICT) for continuous fictitious play if every pair of points in $L$ can be connected by finitely many arbitrarily long pieces of orbits of CFP lying completely within $L$ with arbitrarily small jumps between them.\footnote{For the formal definitions of {\it attractor} and {\it attractor-free set} see \citet[p.~675]{Benaim06}; for the definition of {\it ICT} see \citet[p.~337]{Benaim05}. Note that the definition of invariance in \citet{Benaim05, Benaim06} applies to the best-reply dynamics, so an appropriate time rescaling must be used to apply it to CFP (see footnote \ref{ft:BRD}). This explains that their definition considers solutions defined for all $t \in \R$ while ours considers solutions defined for all $t >0$.}  Every ICT set is invariant and attractor free \citep[Property 2]{Benaim06}. 
The \emph{limit set of the beliefs} of a trajectory $z(t)$ on $\Delta(A_1 \times A_2)$ is the set of all accumulation points of its marginals $(z_1(t), z_2(t)) \in \Delta(A_1) \times \Delta(A_2)$ as $t \to  \infty$. 

\begin{thm}  \label{thm:main} For every no-regret dynamics in class $\mathcal R$, the limit set of the beliefs is almost surely internally chain transitive for continuous fictitious play.\footnote{In the statement of Theorem \ref{thm:main}, CFP can be replaced by the best-reply dynamics since they clearly have the same ICT sets (see footnote \ref{ft:BRD}).} 
\end{thm}
We give here a sketch of the proof. The details are given in Appendix A.1. A discrete-time trajectory $(x_1(t),x_2(t))_{t=1}^{\infty}$ on $\Delta(A_1)\times\Delta(A_2)$ is a \emph{payoff perturbed} DFP trajectory if there exists a positive sequence $(\eps_t)$  converging to zero such that (\ref{Eq:FP}) holds and $a_{i}(t)$ is a vertex of $\Delta(A_i)$ associated with a pure $\eps_t$-best reply to $x_{-i}(t-1)$, for all $i=1,2$ and all $t>1$. A no-regret dynamics in class $\mathcal{R}$ generates a trajectory $(z(t))_{t=1}^{\infty}$ on $\Delta(A)$ and an associated \emph{sequence of beliefs} $(z_1(t), z_2(t))$ on $\Delta(A_1) \times \Delta(A_2)$.  Building on Lemma \ref{lm0}, we show that this sequence of beliefs is almost surely a payoff perturbed DFP trajectory. By an auxiliary lemma, this implies that this is almost surely a \emph{graph-perturbed} DFP trajectory: a notion similar to payoff-perturbed trajectory, but for another definition of perturbed best-reply, the one used in the theory of perturbed differential inclusions \citep{Benaim05, Benaim06}. It follows that the continuous-time interpolation of this sequence of beliefs is almost surely a perturbed solution of CFP, in the sense of \citet{Benaim05}. Theorem \ref{thm:main} then follows from Theorem 3.6 of \citet{Benaim05}.

Since ICT sets are invariant, a consequence of Theorem \ref{thm:main} is the following:
\begin{cor}
\label{cor:main}
Let $\cal{A}$ be the global attractor of CFP  (i.e., its maximal invariant set, see \citet{Benaim05}). For any no-regret dynamics in class $\mathcal R$, the limit set of the beliefs is almost surely a subset of $\cal{A}$.  
\end{cor}
Note the similarity with Propositions 5.1 and 5.2 of \citet{Hofbauer09}, who study the links between the time-average of the replicator dynamics and CFP.

\subsection{Applications of Theorem \ref{thm:main} and comments.}\label{s:cons}

Theorem \ref{thm:main} allows for alternative and sometimes much shorter proofs of most known convergence properties of no-regret dynamics.
Below, we write that no-regret dynamics converge to some set $E$ if the limit set of the beliefs is almost surely a subset of $E$.\footnote{Note that some applications of Theorem \ref{thm:main} (points (a), (b) and (c) below) lead to the same conclusions about no-regret dynamics as those about the time average of the replicator dynamics described in \citet[p.~267, points (2), (3) and (4)]{Hofbauer09}.}

\medskip\noindent {\bf (a)} For any game which is best-reply equivalent to a two-person zero sum game, 
the global attractor of CFP is the set of Nash equilibria \citep{Hofbauer06}. Hence all no-regret dynamics 
in class $\mathcal R$ converge to the set of Nash equilibria. Actually, in zero-sum games, if the correlated 
action $z$ is in the Hannan set (recall that this is the set of correlated actions that satisfy {\it no-regret} for all 
players), then $(z_1, z_2)$ is a Nash equilibrium. Consequently, in zero-sum games 
all dynamics that lead to no regret (not only those in class $\mathcal R$) converge to the set of Nash equilibria. 
This holds more generally for {\it stable bimatrix games} \citep{Hofbauer09JET}, because these are rescaled zero-sum games in the sense of \citet{Hofbauer98}, as is easily shown and was known to Josef Hofbauer (private communication). 

\medskip\noindent {\bf (b)}  For games with strictly dominated strategies, the global attractor of CFP is contained in the face of the simplex with no weight on these strategies. Hence all no-regret dynamics in class $\mathcal R$ converge to this face. Similarly, these dynamics converge to the unique Nash equilibrium in strictly dominance solvable games.

\begin{figure}[!htb]\centering
\begin{tabular}{cc}
$
\begin{array}{c|ccc}
& A & B & C\\
\hline
A & 2 & 1 & -4 \\
B & 1 & 0 & -1 \\
C & -4 & -1 & -2
\end{array}
$ & \qquad \qquad
$
\begin{array}{l|llll}
& A & A^{-} & B & B^{-} \\
\hline
A   & 1   & 1   & 0   & 0\\
A^{-} & 1-\eps &1-\eps  & -\eps  & -\eps \\
B   & 0  & 0     & 1  & 1\\
B^{-} & -\eps&  -\eps   & 1-\eps& 1-\eps\\
\end{array}
$ \\
(i) &  \qquad \qquad (ii) 
\end{tabular}
\caption{}\label{F:H}
\end{figure}

Contrary to (a), this need not be true for all dynamics that lead to no regret. Indeed, convergence to the Hannan set or even to the reduced Hannan set does not guarantee elimination of strictly dominated strategies. Consider, for instance, the games shown on Fig.~\ref{F:H}. Both games are symmetric, so we indicate only the payoffs of the row player. Game (i) is an identical interest game which is strictly dominance solvable;  yet the correlated action putting probabilities $1/3$ on each diagonal square is in the reduced Hannan set. For $\eps=0$, game (ii) is a coordination game with duplicate strategies. For $\eps>0$, the duplicates $A^{-}$, $B^{-}$ are penalized and become strictly dominated. Thus, the correlated action putting probability $1/2$ on $(A^{-}, A^{-})$ and $1/2$ on $(B^{-}, B^{-})$ puts only weight on strictly dominated actions. Yet, for $\eps \leq 1/2$, it belongs to the Hannan set.\footnote{See also the game of \citet[p.~205]{Moulin78}, where the third strategy of player 1 is strictly dominated but has a positive marginal probability under some correlated actions in the Hannan set.}

\medskip\noindent {\bf (c)} In weighted potential games, all internally chain transitive sets of CFP are (subsets of) connected components of Nash equilibria on which the payoffs are constant \citep[see][Theorem 5.5 and Remark 5.6]{Benaim05}. Hence by Theorem \ref{thm:main}, all no-regret dynamics in class $\mathcal R$ converge to such components. Note that the original proof is much longer \citep[Appendix A]{Hart03}.

\medskip\noindent {\bf (d)} If the beliefs $(z_1(t),z_2(t))$ of a no-regret dynamics converge to the set of Nash equilibria, then the average realized payoff converges to the set of Nash equilibrium payoffs. To see why this is true, let $\hat{z} \in\Delta(A)$ be a limit point of $\{z(t)\}$ and let the marginals $(\hat{z}_1,\hat{z}_2)\in\Delta(A_1)\times\Delta(A_2)$ constitute a Nash equilibrium. By Corollary \ref{Cor:1} the maximal regret converges to zero, so for every $i=1,2$
\begin{align*}
u_{i}(\hat z) &= \max_{k\in A_i} u_{i}(k, \hat{z}_{-i})= u_{i}(\hat{z}_i, \hat{z}_{-i}).
\end{align*}
This result illuminates an important difference between no-regret dynamics and discrete fictitious play. It is well known that under DFP, if the beliefs of the players converge to a Nash equilibrium, their average realized payoffs \emph{need not} approach the set of Nash equilibrium payoffs,  whereas under no-regret dynamics it is always the case.

\begin{figure}[!htb]
\[
\begin{array}{c|ccc}
& A & B & C\\
\hline
A & 0,0 & 1,0 & 0,1 \\
B & 0,1 & 0,0 & 1,0 \\
C & 1,0 & 0,1 & 0,0
\end{array}
\] 
\caption{}\label{F:S} 
\end{figure}

\medskip\noindent {\bf (e)} Consider the $3\times 3$ game of Fig.~\ref{F:S} due to \citet{Shapley64}, the historical counterexample to the convergence of fictitious play. This game has a unique equilibrium, in which both players randomize uniformly. Though this equilibrium attracts some solutions of continuous fictitious play (e.g. all those that start and remain symmetric), almost all solutions converge to a hexagon, the so-called Shapley polygon \citep{Gaunersdorfer95, Shapley64, Sparrow08}. It may be shown that the only ICT sets are the Nash equilibrium and the Shapley polygon. Consequently, the limit set of any no-regret dynamics in class $\mathcal R$ is almost surely one of these two sets.

\medskip\noindent {\bf (f)} In a number of classes of games, convergence of discrete fictitious play to the set of Nash equilibria has been established, but analogous results for continuous fictitious play are lacking. 
Thus we cannot use Theorem \ref{thm:main}. These classes of games include generic $2\times n$ games \citep{Berger05}, generic ordinal potential games, quasi-supermodular games%
\footnote{Also known as games of strategic complementarities (e.g., \citet{Tirole88}).} with diminishing returns \citep{Berger07}, and some other special
classes (see, e.g., \citet[p.~260]{Sparrow08}). For ordinal potential games and quasi-supermodular games with diminishing returns, \citet{Berger07} proves convergence to the set of Nash equilibria of \emph{some} solutions of continuous fictitious play  as defined by (\ref{Eq:CFP2}) (see our footnote \ref{ft:BRD}). 
This is not enough to apply the results of \citet{Benaim05}. The same problem arises in \citet{Krishna98}. Actually, as explained below, convergence of CFP to the set of Nash equilibria would not suffice to use Theorem \ref{thm:main}: we would need some additional structure, such as a Lyapunov function, to get more information on the ICT sets. 

\medskip\noindent {\bf (g)} 
Consider a bimatrix game in which all solutions of CFP converge to the set of Nash equilibria. Because the definition of attractor requires uniform attraction, this does not imply that the set of Nash equilibria is an attractor. Neither does it imply that all ICT sets are contained in the set of Nash equilibria, as shown in Appendix A.2. Therefore, we cannot deduce from Theorem \ref{thm:main} that no-regret dynamics in class $\mathcal{R}$ converge to the set of Nash equilibria; whether this is always the case remains an open question.

\medskip\noindent {\bf (h)} We show in Section \ref{sec:cont} that Theorem \ref{thm:main} also applies, and under weaker assumptions, to the continuous-time version and to the expected version of no-regret dynamics in class $\mathcal{R}$. As apparent from the proof, the existence of a potential is not essential: for a good choice of behavior when there are no regrets,  Theorem \ref{thm:main} holds for any no-regret dynamics such that a player always chooses an action with positive regret whenever he has one. It also applies to certain no-regret dynamics that do not have this property, such as the exponential weight algorithm (see Remark \ref{rem:ew} at the end of Appendix A.1).

\medskip\noindent {\bf (i)}
Let us now comment on extensions of our results to $n$-player games. 
The definition of no-regret dynamics, as well as Proposition \ref{Thm:HM}, extend to the $n$-player setting straightforwardly (e.g., \citet{Hart01Gen}).
The appropriate extension of CFP is {\it correlated CFP} where at each time $t$ every player chooses a best reply action to the {\it correlated} past average play of the others. Specifically, an absolutely continuous function $z:[1,\infty)\to\Delta(A)$ is a solution of correlated CFP if it is almost everywhere differentiable and satisfies
\[
\dot z(t)\in \frac{1}{t}(\overline{BR}(z(t))-z(t)),
\]
where the correlated best-reply correspondence $\overline{BR}: \Delta(A) \rightrightarrows \Delta(A)$ is defined by $BR(z)=\times_{i=1}^n \overline{BR}_i(z_{-i})$ where $\overline{BR}_i(z_{-i})$ is the set of mixed best replies of player $i$ to the correlated action $z_{-i}\in \Delta\left(A_{-i} \right)$. 

In $n$-player games with linear incentives \citep{Selten95}, also known as polymatrix games \citep{Yanovskaya68},  the correlated and independent best-reply correspondences coincide; that is, for any correlated action $z \in \Delta(A)$, $\overline{BR}(z)=BR((z_1,..,z_n))$ where $(z_1,..., z_n)$ is the vector of marginals of $z$, and $BR$ the standard (independent) best-reply correspondence. For such games, Theorem \ref{thm:main} extends easily. However, this is not the case in general. The main problem is that the correlated best reply correspondence is not convex valued; that is, $\overline{BR}(z)$ is not in general a convex subset of $\Delta(A)$.\footnote{This is due to the fact that elements of $\overline{BR}(z)$ are independent distributions and that the average of two independent distributions need not be an independent distribution.} This creates two issues: 

(i) Existence of solutions of correlated CFP is not guaranteed by the classical results on differential inclusions we are aware of (e.g., \citet{Aubin84}).

(ii) The theory of perturbed differential inclusions \citep{Benaim05} does not apply to non-convex valued correspondences. 

The first issue can be solved by building piecewise linear solutions of correlated CFP following the same ideas as for two-player games (see \citet{Hofbauer95}).\footnote{Assuming that $z(t)$ is well defined, call $G_r(t)$ the game in which the players are reduced to their 
best-replies to $z(t)$. Start with some initial condition $z(t_0)$. Then point to a Nash equilibrium 
of $G_r(t_0)$ (i.e. fix $b \in NE(G_r(T_0))$ and choose $q(t) =b$) till the first time, $t_1$, when, for some 
player $i$,  a strategy which was not a best-reply to $z(t_0)$ is a best-reply to $z(t_1)$. Then iterate. 
If the times $t_n$ accumulate towards some time $t^{\ast}$, then use the fact that $z(t)$ must have a 
limit when $t \to t^{\ast}$ (because $z(t)$ is Lipschitz). Call it $z(t^{\ast})$ and restart from $z(t^{\ast})$. Note that there might in principle be a countable infinity of such accumulation points $t^{\ast}$, and 
that they might themselves accumulate in some point $t^{\ast \ast}$, but then define $z^{\ast \ast}$ 
as before and restart from there, etc. The largest (forward time) interval on which such a solution can be built is both 
open and closed in $[t_0, +\infty)$ and is thus equal to $[t_0, +\infty)$.} Moreover, due to Remark \ref{rem:uni}, Proposition \ref{P:CFP} 
extends to the $n$-player setting. It then asserts that correlated CFP leads to no regrets. Lemma \ref{lm0} also extends: 
it asserts that if the maximal regret of player $i$ is less than $\eps$, then she plays only $\eps$-best reply actions to the {\it correlated} average play of the opponents. It follows that, analogously to two-player games, interpolated trajectories of no-regret dynamics are almost surely perturbed solutions of correlated CFP.
However, we cannot proceed to an analog of Theorem \ref{thm:main} because of the second issue. Thus, whether there is a formal relation between no-regret dynamics and correlated CFP in $n$-player games remains an open question. Similarly, the results of \citet{Hofbauer09} on the links between the time-average of the replicator dynamics and CFP are restricted to bimatrix games (or games with linear incentives).

\section{Curb sets}\label{Sec:curb}

Theorem \ref{thm:main} does not answer whether {\it attracting sets} of CFP have an analogous property under no-regret dynamics. 

A set $\mathcal C\subset \Delta(A)$ is {\it eventually attracting} under a no-regret dynamic process if with any given probability it captures all no-regret trajectories originating from a small enough neighborhood of $\mathcal C$ at all distant enough periods. Formally, $\mathcal C$ is eventually attracting if for every $\pi>0$ there exists $\eps>0$ and a period $T$ such that: 
for every $t_0 \geq T$, if $z(t_0)$ is in an $\eps$-neighborhood of $\mathcal C$, then $z(t)$ converges to set $\mathcal C$ with probability at least $1-\pi$.\footnote{We say that $z(t)$ converges to $\mathcal C$ if $\inf_{c\in \mathcal C} ||z(t)-c||\to 0$ as $t\to\infty$.}

For this section it is convenient to replace assumption (Q2) by the following one:
\begin{equation}
\begin{array}{l}
\text{If a player's maximal regret is nonpositive, then she plays a best-reply} \\ \text{to the empirical distribution of her opponent.} 
\end{array}\tag{Q2$'$}
\end{equation}
This is not essential, since the interesting histories are those where both players have positive regrets, in which case (Q2) plays no role.\footnote{Recall that by Remark \ref{rem:Q1}, if a player has positive maximal regret, then it remains positive forever. So we can consider histories from a distant enough period $t_0$ where both players have positive regrets and (Q2) plays no role. If $t_0$ does not exist, i.e., some player {\it always} has nonpositive maximal regret, then Proposition \ref{Thm:HM} and (Q2) imply that her play is constant, whereas her opponent's play must approach a best reply to it, leading to Nash equilibrium. By replacing (Q2) by (Q2$'$) we avoid dealing with this issue.}

A strict Nash equilibrium is eventually attracting. Indeed, if $z(t_0)$ is close enough to a vertex of $\Delta(A)$ corresponding to a strict Nash equilibrium $a=(a_1,a_2)$,
then for each player $i$, action $a_i$ is the unique best reply and there is a negative regret for any action other than $a_i$. 
Since by (R3) only actions with positive regret can be chosen, and by (Q2$'$) only best-reply actions can be chosen if {\it all} regrets are nonpositive, action $a_i$ will be played by each player $i$ in the following period, and so on.

Let us now consider a standard generalization of strict Nash equilibria.  
For each $i=1,2$, let $B_i\subset A_i$. With a slight abuse of notation, 
denote by $\Delta(B_i)$ the set of probability measures on $A_i$ with support on $B_i$ only. 
The product set $B=B_1\times B_2$ is \emph{closed under rational behavior (curb)} (\citet{Basu91}) if
\[BR_i(x_{-i})\subset\Delta(B_i) \ \text{whenever $x_{-i}\in \Delta(B_{-i})$},\  i=1,2.
\]
That is, the set $B$ is curb if the players' pure best reply profiles are contained in $B$ whenever they believe that no actions
outside of $B$ should be played. 

Curb sets are known to be attracting under CFP (e.g., \citet[Lemma 7]{Balkenborg12}). However, they need not be attracting under no-regret dynamics in class $\mathcal{R}$. Indeed, even if the support of $z(t_0)$ is contained in some curb set $B$, there may be positive regrets for actions outside of $B$, since $B$ need not be closed under \emph{better} replies. 
However, we show that the intersection of the Hannan set and the set of correlated actions with support on a curb set is eventually attracting. 

Formally, let $B=B_1\times B_2$ be a curb set. Let $\Delta_B(A)$ denote the set of correlated actions with support on $B$ only. Let $H_B=H\cap\Delta_B(A)$.

\begin{prop}\label{Prop:Curb}
For every curb set $B$, the set $H_B$ is eventually attracting under every no-regret dynamics in $\mathcal R$.
\end{prop}
The proof is based on the following observations. For every curb set $B$, if the average play is close enough to $H_B$, then regrets for all actions outside of $B$ are negative (since $B$ is curb). Hence, by condition (R3), only actions in $B$ will be played in the immediate future. On the other hand, almost sure convergence of maximum regret to zero suggests that, so long as the players choose only actions in $B$, the average play will approach $H_B$, thus reinforcing the former observation. To prove the result, however, we need to establish bounds on the maximal future regret {\it conditional on certain histories} (namely, conditional on being close to  $H_B$) that \citet{Hart01Gen} do not provide. The complete proof is relegated to Appendix A.3.

\section{Continuous-time and expected no-regret dynamics} \label{sec:cont} 

We now prove an analog of Theorem \ref{thm:main} for continuous-time dynamics \citep{Hart03} and the expected version of discrete-time dynamics. Both describe trajectories of average {\it intended (mixed)} play, rather than average realized (pure) play. For this reason, 
condition (R3) is not needed. Indeed, the interest of (R3) is that, together with (Q1), it 
requires every realized action to be a better reply to the opponents empirical distribution 
of play (whenever such actions exist). But now we only need every mixed (expected) action to be 
a better reply, and this follows already from conditions (R1)-(R2) and (Q1). Besides, these dynamics are deterministic, hence the results we obtain hold {\it surely} (not just {\it almost surely}). The proofs are based on Appendix A.1 and are best understood after reading it.

Consider a continuous-time dynamics 
\begin{equation}
\label{eq:contz}
\dot z(t)=\frac{1}{t}(\bar{q}(t)-z(t))
\end{equation}
 where $\bar{q}(t)= q_{1}(t)\otimes q_2(t) \in \Delta(A)$ is the (independent) joint play at time $t$ and $z(t)$ the average correlated play. There are two differences with (\ref{Eq:DT}): time is now continuous, and, more importantly, realized play $a(t)$ has been replaced by intended mixed play $\bar{q}(t)$. As in CFP, start at time $1$ with some initial condition $z(1) \in \Delta(A)$. Assume that whenever $R_{i,\max}(t) >0$
\begin{equation}\label{eq:contq} 
q_{i,k}(t)=\frac{\nabla_k P_{i}(R_{i}(t))}{\sum_{s\in A_{i}}\nabla_{s} P_{i}(R_{i}(t))},\text{ \ }k\in A_{i}
 \end{equation} 
 where $P_i$ is a $C^1$ potential function satisfying (R1), (R2) and the technical condition:\bigskip\\ 
 \noindent 
\textbf{(P4$'$)} \quad There exists $0<\rho_2<\infty$ such that $\nabla P_i(x) \cdot x \leq \rho_2 P_i(x)$ for all $x \in  \R^{A_i}_{-}$.\bigskip\\ 
This is a part of condition (P4) in \citet{Hart03}. 
\begin{prop} \label{prop:cont}
 Let $z(t)$ be a solution of (\ref{eq:contz}) and (\ref{eq:contq}) with $P_i$ satisfying conditions \emph{(R1)}, \emph{(R2)} and \emph{(P4$'$)} for all $i=1,2$. Assume that the initial condition $z(1)$ is such that both players have some positive regrets: $R_{i,\max}(1) >0$ for all $i=1,2$. Then the limit set of the beliefs is internally chain transitive for continuous fictitious play.
 \end{prop}
 
\begin{pf} Let $\eps_i(t):=R_{i,\max}(t)$. \citet[Theorem 3.1 and Lemma 3.3\footnote{Note a typo in the proof of Lemma 3.3 in \citet{Hart03}: (P3) should be replaced by (P4). Moreover, only our condition (P4$'$) is used in the proof of Lemma 3.3 in \cite{Hart03}.}]{Hart03}  
show that if $\eps_i(1)>0$, then $\eps_i(t)>0$ for all $t$, and $\eps_i(t) \to 0$ as $t \to +\infty$. Moreover, by (R2) applied to $x=R_i(t)$ and definition of $q_i$, we have: $u_i(q_i, z_{-i}) - u_i(z) = q_i \cdot R_i >0$ (this is equation (3.3) in \citep{Hart01Gen}).
Thus by Lemma \ref{lm0}, $q_i \in BR_i^{\eps_i(t)}(z_{-i})$. Together with Lemma \ref{lm2} in Appendix A.1, this implies that $(z_1(\cdot), z_2(\cdot))$ is a perturbed solution of CFP in the sense of \citet{Benaim05}. The result then follows from Theorem 3.6 of \citet{Benaim05}.     
\end{pf}

\begin{rem} Assume that if all initial regrets of a player are nonpositive then the dynamics is defined as in \citet{Hart03}, equation (4.9). Then it is easily seen that the result of Proposition \ref{prop:cont} holds for any initial condition $z(1)$.
\end{rem}

\noindent {\bf Expected discrete-time dynamics. } The expected motion in (\ref{Eq:DT}) is described by 
\begin{equation*}
z(t)-z(t-1)=\frac{1}{t}\left(\bar{q}(t)-z(t-1)\right).
\end{equation*}
where $\bar{q}(t)=q_1(t)\otimes q_2(t)$ is the expectation of $a(t)$. 
Assume that $q_i$ is derived by (Q1) from a potential function satisfying (R1)--(R2). Let $\eps_i(t):=R_{i,\max}(t)$. It is easily seen that, as for continuous-time dynamics, $\eps_i(t) \to 0$ as $t \to +\infty$, and if $\eps_i(1)>0$, then for all $t$, $\eps_i(t)>0$ and $q_i \in BR_i^{\eps_i(t)}(z_{-i})$. Due to Lemmata \ref{lm2} to \ref{lm4} of Appendix A.1 and to Theorem 3.6 of \citet{Benaim05}, it follows that for a good choice of behavior when all regrets are initially nonpositive, the limit set of the beliefs is internally chain transitive for CFP.

\section*{Appendix}
\renewcommand{\thesubsection}{A.\arabic{subsection}}

\subsection{Proof of Theorem \ref{thm:main}} \label{app:thm}

Denote by $\widehat{BR}{}_i^{\eps}(x)$ the correspondence whose graph is the $\eps$-neighborhood of the graph of $BR_i$:
\begin{equation*}
\widehat{BR}{}_i^{\eps} (x_{-i}) = \left\{ x_i \in \Delta(A_i)  \ \left| \ \begin{array}{l}
\exists (x_i^{\ast}, x_{-i}^{\ast}) \in \Delta(A_1) \times \Delta(A_2) \  \text{s.t.} \\
x_i^{\ast} \in BR_i(x_{-i}^{\ast}), \ \text{and} \
||(x_i^{\ast}, x_{-i}^{\ast}) - (x_i,x_{-i})||_{\infty} \leq \eps
 \end{array}
  \right.\right\}
\end{equation*}
Let $\widehat{BR}{}^{\eps}(x)=\widehat{BR}{}_1^{\eps} (x_{2})\times \widehat{BR}{}_2^{\eps} (x_{1})$. In words, action $x_i$ is an {\it $\eps$-graph perturbed best reply to $x_{-i}$} if there is an action $\eps$-close to $x_i$ which is an exact best-reply to an action $\eps$-close to $x_{-i}$. This is the notion of perturbation used in the theory of perturbed differential inclusions (\citet{Benaim05, Benaim06}). As illustrated by the example below, 
it is different from the notion of  perturbation of payoffs in the $\eps$-best reply correspondence, i.e. $BR^{\eps}(x)=BR_1^{\eps}(x_2) \times BR_2^{\eps}(x_1)$ with 
\begin{equation*}
BR_i^{\eps} (x_{-i})= 
\Big\{x_i \in \Delta(A_i) \ \Big| \ u_i(x_i, x_{- i})  \geq \max_{k\in A_i} u_i(k,x_{-i})-\eps \Big\}, \ \ i=1,2.
\end{equation*}

\begin{figure}[!htb]
\[
\begin{array}{c|cc}
& L & R \\
\hline
T & 1  & 0 \\
C & 0 & 1 \\
B & \frac{1}{2}-\eta & \frac{1}{2}-\eta
\end{array}
\]
\caption{}
\label{F:EMP2}
\end{figure}
Consider a game where the payoffs of player 1 are given by Fig~\ref{F:EMP2}.
Let $\eps \in(0,1/2)$ and let $x_2^{\eps}=\left(\frac{1}{2} + \eps\right) L + \left(\frac{1}{2} - \eps\right) R$. 
The pure action $C$ is a $2\eps$-best reply to $x_2^{\eps}$. Using the sup norm, it is at distance 1 from pure action $T$, the unique exact best reply to $x_2^{\eps}$. Nevertheless, $C$ is an $\eps$-graph perturbed best reply, because it is an exact best reply to $x_2^0$, which is $\eps$-close (in sup norm) to $x_2^{\eps}$. 
By contrast, for all $\eta>0$, action $B$ is an $(\eps+\eta)$-best reply, but only a $1$-graph perturbed best reply to $x_2^{\eps}$. 

A discrete-time trajectory $(x_1(t),x_2(t))_{t=1}^{\infty}$ on $\Delta(A_1)\times\Delta(A_2)$ is a \emph{payoff perturbed} fictitious play trajectory if there exists a positive sequence $(\eps_t)$  converging to zero such that 
\begin{equation*}
x(t)-x(t-1)= \frac{1}{t}\left(q(t)-x(t-1)\right)
\end{equation*}
with $q(t)=(q_1(t),q_2(t))$ and $q_i(t)\in BR^{\eps_t}_i(x_{-i}(t-1))$ for all $i=1,2$ and all $t>1$. It is a \emph{graph perturbed} fictitious play trajectory if the same holds but replacing  $BR^{\eps_t}_i$ with $\widehat{BR}{}^{\eps_t}_i$. A trajectory $(z(t))_{t=1}^{\infty}$ on $\Delta(A)$ generates a \emph{sequence of beliefs} $(z_1(t), z_2(t))$ in $\Delta(A_1) \times \Delta(A_2)$ .

The proof goes as follows. Lemma \ref{lm1} shows that the sequence of beliefs generated by a no-regret dynamics is a payoff perturbed FP trajectory. Together with Lemma \ref{lm2}, this implies that it is a graph-perturbed FP trajectory (Lemma \ref{lm3}). It follows that the interpolated process of a no-regret dynamics trajectory is a perturbed solution of CFP (Lemma \ref{lm4}). The result then follows from \citet{Benaim05}.    
\begin{lem}
\label{lm1} The sequence of beliefs of a solution of a no-regret dynamics in class $\mathcal{R}$ is almost surely a payoff perturbed DFP trajectory.  
\end{lem}
\begin{pf}  If $R_{i,\max}(t) \leq 0$ for all $t$, then by Remark \ref{rem:Q2}, player $i$ plays an $\eps(t)$-best reply for some $\eps(t)$ converging to zero. Otherwise, $R_{i,\max}(t_0)>0$ for some $t_0 \in \N^{\ast}$. Then for all times $t >t_0$, $R_{i,\max}(t)>0$ (by Remark \ref{rem:Q1}) and player $i$ plays an $R_{i,\max}(t)$-best reply  by Lemma \ref{lm0} and conditions (R3) and (Q1). Since $R_{i,\max}(t) \to 0$ almost surely, the result follows.   
\end{pf}

\begin{lem} \label{lm2}
 Let $X$ be a compact subset of $\R^m$ and $F$ a correspondence from $X$ to itself. For any $\delta \geq 0$, let $\hat{F}_{\delta} : X \rightrightarrows X$ denote the correspondence whose graph is the $\delta$-neighborhood of the graph of $F$:
\[ 
\hat{F}_{\delta} (x) = \Big\{y \in X\,\Big|\,\exists (x^{\ast}, y^{\ast})\in X^2\ \text{s.t.}\  y^{\ast} \in F(x^{\ast}) \ \text{and} \ ||(x^{\ast}, y^{\ast}) - (x,y)||_{\infty} \leq \delta \Big\}. 
\] 
For any $\alpha>0$, let $G_{\alpha}$ be an u.s.c. correspondence from $X$ to itself. Assume that for each $x$ in $X$: 
 
 (i) $\alpha < \alpha' \Rightarrow G_{\alpha}(x) \subset G_{\alpha'}(x)$ (that is,  $(G_{\alpha})_{\alpha >0}$ is increasing w.r.t. inclusion);
 
 (ii) $\bigcap_{\alpha >0} G_{\alpha}(x) \subset F(x)$.
 
\noindent Then for every $\delta>0$ there exists $\alpha>0$ such that for each $x$ in $X$,  $G_{\alpha}(x) \subset \hat{F}_{\delta} (x)$. 
\end{lem}
\begin{pf}  By contradiction, assume that there exists $\delta>0$, a decreasing sequence $(\alpha_n)$ converging to zero, and sequences  
$(x_n)$ and $(y_n)$ of points in $X$ such that $y_n \in G_{\alpha_n}(x_n)\backslash \hat{F}_{\delta}(x_n)$ for all $n$. 
By compactness of $X$, we can assume that $(x_n)$ and $(y_n)$ converge respectively to $x^{\ast}$ and $y^{\ast}$. Fix $k \in \N$. For all $n \geq k$, $y_n \in G_{\alpha_n}(x_n) \subset G_{\alpha_k}(x_n)$ by (i). Since $G_{\alpha_k}$ is u.s.c., it follows that $y^{\ast} \in G_{\alpha_k}(x^{\ast})$. Therefore, by (i) and (ii)
\[y^{\ast} \in \bigcap_{k \in \N} G_{\alpha_k}(x^*) = \bigcap_{\alpha >0} G_{\alpha}(x^*) \subset F(x^*)\] 
But for $n$ large enough, $||(x^{\ast}, y^{\ast}) - (x_n, y_n)||_{\infty} < \delta$, hence $y_n \in \hat{F}_{\delta} (x_n)$, a contradiction.  
\end{pf}

\medskip
Applied to the best-reply correspondence, Lemma \ref{lm2} implies that for any $\delta>0$, 
an $\eps$-perturbed best-reply is a $\delta$-graph perturbed best-reply, provided $\eps$ is small enough. Thus we have the next result.

\begin{lem} \label{lm3}  Any payoff perturbed DFP trajectory is a graph perturbed DFP trajectory.
\end{lem}  
\begin{pf} Let $\eps_t \to 0$. Let 
\[
\delta_t=\min \left\{ \delta \geq 0 \,\left| \, \forall i =1,2,\, \forall x \in \Delta(A_1) \times \Delta(A_2), BR_i^{\eps_t} (x_{-i}) \subset  \widehat{BR}{}_i^{\delta}(x_{-i})\right.\right\}.
\]
Applying Lemma \ref{lm2} with $X=\Delta(A_1) \times \Delta(A_2)$, $G_{\eps}=BR^{\eps}$ and $F=BR$, we obtain that $\delta_t \to 0$. The result follows.
\end{pf} 

\medskip
Given a discrete-time trajectory $x(n)=(x_1(n), x_2(n))$ on $\Delta(A_1) \times \Delta(A_2)$, with $n\in\N^{\ast}$, define its interpolated process $x : [1, +\infty) \to \Delta(A_1) \times \Delta(A_2)$ as follows. For all $t \in [n, n+1)$ let 
$tx(t)=nx(n) + (t-n)q(n)$, where $q_i(n)=(n+1)x_i(n+1)- nx_i(n)$, $i=1,2$. This is equivalent to
 \[x_i(t)- x_i(n)=\frac{t-n}{t} (q_i(n)-x_i(t)), \ \ i=1,2.\]
Hence for all $t \in (n, n+1)$ we have $||x(t) - x(n) ||_{\infty} \leq \frac{1}{n+1}$ and 
\begin{equation}
\label{eq:int}
\dot{x}(t)= \frac{1}{t}(q(n) - x(t))
\end{equation}
An absolutely continuous function $x: [1, +\infty) \to \Delta(A_1) \times \Delta(A_2)$ is a perturbed solution of CFP if there exists a vanishing function $\eps : \R_+ \to \R_+$ such that for almost all $t$, 
\begin{equation}\label{eq:ptbCFP} 
\dot{x} \in \frac{1}{t} \left(\widehat{BR}{}^{\eps(t)} (x)- x\right) \quad \mbox{ where $x=x(t)$.}
\end{equation}
\begin{lem} \label{lm4} 
The interpolated process of a graph perturbed DFP trajectory is a perturbed solution of CFP.
\end{lem}
\begin{pf}
Consider a discrete time trajectory $(x_1(n),x_2(n))_{n\in\N}$ such that
\[
x_i(n)-x_i(n-1)= \frac{1}{n}\left(q_i(n)-x_i(n)\right), \ \ i=1,2,
\]
with $q_i(n)\in \widehat{BR}{}_i^{\eps_n}(x_{-i}(n-1))$ and $\eps_n \to 0$. For all $n$ and all $t \in [n, n+1)$, let $\eps(t)=\eps_n + 2/n$. 
Obviously, $\eps(t) \to 0$ as $t \to \infty$. Moreover, for all $t \in (n, n+1)$, the interpolated process satisfies $||x_{-i}(t) - x_{-i}(n-1) ||_{\infty} \leq \frac{1}{n+1} + \frac{1}{n} < 2/n$, 
so $q_i(n) \in \widehat{BR}{}^{\eps(t)}_i(x_{-i}(t))$. Therefore (\ref{eq:int}) implies (\ref{eq:ptbCFP}) (see also \citet[Proposition 2.2]{Faure10}). 
\end{pf}

\medskip
We can now prove Theorem \ref{thm:main}. By Lemmata \ref{lm1} and \ref{lm3}, the sequence of beliefs of a solution of a no-regret dynamics in class $\mathcal{R}$ is almost surely a graph perturbed DFP trajectory. Hence, by Lemma \ref{lm4}, its interpolated process $x(t)$ is a perturbed solution of CFP. This implies that $x(e^{t})$ is almost surely a perturbed solution of the best-reply dynamics, in the sense of \citet[Definition II]{Benaim05}. Theorem \ref{thm:main} now follows from Theorem 3.6 of \citet{Benaim05}.\footnote{The definition of perturbed solution in \citet{Benaim05} is different from ours but equivalent.}

\begin{rem}
\label{rem:ew}
Assume that at stage $t$, for each $i=1,2$, player $i$ chooses a pure action according to a mixed action $q_i(t)$  that depends on the previous history $h(t-1)$. Do not assume conditions (R1)--(R3) and (Q1), but assume that there exists a vanishing sequence $(\eps_t)$ such that for for all $t>1$ and any previous history $h(t-1)$, $q_i(t) \in BR_i^{\eps_t}(z_{-i}(t-1))$, $i=1,2$.  Then it follows from Lemma \ref{lm2}, the above proof and \citet[Proposition 1.4 and a variant of Proposition 1.3]{Benaim05} that Theorem \ref{thm:main} applies. 
As is well known, this is the case for  the exponential weights algorithm \citep{Freund99,Littlestone94} that corresponds to 
\[q_{i,k}(t):=\frac{\exp{\beta_t u_i(k,z_{-i})}}{\sum_{s \in A_i} \exp(\beta_t u_i(s, z_{-i}))} \]
with $z_{-i}=z_{-i}(t-1)$, $\beta_t \to +\infty$ as $t \to \infty$, and $\beta_t < t^\alpha$ for some $\alpha \in (0,1)$ to ensure that this is a no-regret dynamics (see, e.g., \citet{Benaim11}). The above assumptions are not (or not trivially) satisfied by no-regret dynamics in class $\mathcal{R}$. Indeed, the rate at which the maximal regret vanishes, hence the value $\eps_t$ such that $q_i(t) \in BR_i^{\eps_t}(z_{-i}(t))$, may depend on the trajectory.
\end{rem}

\subsection{ICT sets when all solutions converge to Nash equilibria}
\label{sec:ICT}
The fact that all solutions of the best-reply dynamics converge to the set of Nash equilibria 
does not guarantee that ICT sets contain only Nash equilibria. We provide counterexamples below.

\begin{ex} [single-population dynamics]
\label{ICT:ex1}
Consider the following symmetric $3 \times 3$ game: 
\[
\begin{array}{c|ccc} 
& A & B & C \\
\hline  
A & 0 & 0 & 0 \\
B & 0 & 0 & 0 \\
C & -1 & 0 & 0
\end{array}
\]

Denote a mixed action by $x=(x_A, x_B, x_C)$. Then $(x,x)$ is a Nash equilibrium if and only if  $x_A=0$ or $x_C=0$. 
It is easily seen that all solutions of the best-reply dynamics converge to the set of  symmetric Nash equilibria. However, the whole state space is ICT. Indeed, any mixed action $x$ can be connected to any other mixed 
action $y$ as follows: starting from $x$, follow a solution pointing towards the edge $x_C=0$, then jump on this edge and follow 
a solution pointing towards the pure strategy $B$; once close to $B$, jump on the edge $x_A=0$, and follow a solution pointing towards $C$;  
once close to $C$, make a small jump to reach a point from which a solution points toward $y$; follow this solution and if needed (i.e. if $y_C=0$), 
make one more jump to reach $y$. 

This example is also valid for the replicator dynamics and any payoff monotone dynamics in the sense of, e.g., \citet{Hofbauer96}. The only difference is that traveling from $A$ to $B$ and from $B$ to $C$ cannot 
be done by following solutions of the dynamics but only through long sequences of jumps.  Note also that in an inward cycling Rock-Paper-Scissors game (see e.g., \citet{Hofbauer98}, or \citet{Weibull95}), 
all solutions of the replicator dynamics converge to one of the four rest points but the whole boundary of the state space is ICT (for the replicator dynamics).
\end{ex} 

 \begin{ex}  [$n$-population dynamics]
\label{ICT:ex2}
Similarly, in the bimatrix version of example \ref{ICT:ex1}, all solutions of the two-population best-reply dynamics 
converge to the set of Nash equilibria but the whole state space is ICT. Again, this is true for all payoff monotone 
dynamics. Similar examples can be given for $n$-population dynamics for any 
$n \geq 1$. 

At least for the best-reply dynamics, $2 \times 2$ examples can also be given. Consider, for instance, the $2 \times 2$ game: 
\[
\begin{array}{c|cc} 
& L \quad & R \\  
\hline
T & 0,0 & 0,0  \\
B & 0,0 & -1,0
\end{array}
\]
Denote mixed actions of players $1$ and $2$ by $x=(x_T, x_B)$ and $y=(y_L, y_R)$, respectively.  The set of 
Nash equilibria is the union of the edges $x_B=0$ and $y_R=0$ and all solutions of the two-population best-reply dynamics 
converge to this set. However, the whole triangle $x_T+ y_L \geq 1$ is ICT.  

In these examples, a direct analysis shows that all solutions of no-regret dynamics in class $\mathcal{R}$ converge to the set of Nash equilibria. Thus we do not know whether, in general, convergence of all solutions of CFP to the set of Nash equilibria entails convergence of no-regret dynamics. The point is that this is not guaranteed by Theorem \ref{thm:main}. 
\end{ex}

\subsection{Proof of Proposition \ref{Prop:Curb}}

We need some notation. 
For $z \in \Delta(A)$ and $a \in A$, let $z_a$ denote the probability of $a$ under the correlated action $z$.  
Let $\mathcal U_{\gamma}(H_B)$ be the neighborhood of $H_B$ in which the total weight on action profiles outside of $B$ 
and the potential of each player are below $\gamma$:
\[
\mathcal U_{\gamma}(H_B)= \left\{z \in\Delta(A) \left|
\begin{split}
& \sum\nolimits_{a\notin B} z_a <\gamma, \ \  \text{and}   \\
& P_i(R_i (z)) < \gamma, \ \ i=1,2,
\end{split}
\right.\right\}
\]
where $R_i(z)$ is the regret vector of player $i$:  $R_i(z)=\big(u_i(s,z_{-i})-u_i(z)\big)_{s\in A_i}$.

Let $B=B_1\times B_2$ be a curb set. Let 
\[
\delta_B=\min_{i=1,2}
\min_{z_{-i}\in\Delta(B_{-i})} \left\{\max_{s\in A_i} u_i(s,z_{-i})- \max_{k\in A_i\bs B_i} u_i(k,z_{-i})\right\}
\]
and note that $\delta_B>0$, since $B$ is curb. Now, consider a no-regret dynamics in $\mathcal R$ defined by potentials $P_i$, $i=1,2$, with trajectory $(z(t))_{t \geq 1}$. 
Let $\rho_i(\gamma)$ be the smallest number such that for all $z \in \Delta(A)$
\[
P_i(z) \leq \gamma \ \ \Longrightarrow \ \ 
R_{i, \max} (z) \leq \rho_i(\gamma),
\]
and let $\rho(\gamma)=\max\{\rho_1(\gamma),\rho_2(\gamma)\}$. Let $\gamma_B$ be the solution of
\begin{equation}\label{E:AQ0}
(2 \bar U+\delta_B)\gamma_B+\rho(\gamma_B)-\delta_B=0  
\end{equation}
where $\bar{U}=\max_{i=1,2} \max_{a \in A} |u_i(a)|$ is a payoff bound. 
Since $\rho(\gamma)$ is weakly increasing in $\gamma$ and $\rho(0)=0$, 
there exists a unique solution $\gamma_B$ of \eqref{E:AQ0} and $\gamma_B>0$.

Consider the following event $\mathcal E_{t}$:
\[
   P_i(R_i(t+n))<\gamma_B \ \ \text{for each $i=1,2$ and all $n\in\N^*$}. \tag{$\mathcal E_{t}$}
\]

The statement of Proposition \ref{Prop:Curb} is immediate by the following claims and Corollary \ref{Cor:1}.

\begin{claim} \label{Cm:1}
If $z(t)\in \mathcal U_{\gamma_B}(H_B)$ and event $\mathcal E_t$ holds, then $a(t+n)\in B$ for all $n\in\N^*$.
\end{claim}

\begin{claim} \label{Cm:2}
For every $\pi\in(0,1]$ and every $\gamma\in (0,\gamma_B)$ there exists $t_0$ such that for every $t\ge t_0$, if $z(t)\in \mathcal U_{\gamma}(B)$, then event $\mathcal E_t$ holds with probability at least $1-\pi$.
\end{claim}

For the proof of Claim \ref{Cm:1} we need the following lemma.

\begin{lem}\label{Lem:U}
For any $t$, if $z(t)\in \mathcal U_{\gamma_B}(H_B)$, then $a(t+1)\in B$.
\end{lem}

\begin{pf} Let $z \in \Delta(A)$. If $z$ is close enough to $B$, then 
$\max_{s\in A_i} u_i(s,z_{-i})- \max_{k \in A_i\bs B_i} u_i(k,z_{-i}) > \delta_B/2$ for all $i=1,2$.  
If $z$ is close enough to $H$, then $\max_{s\in A_i} u_i(s, z _{-i})-u_i(z) < \delta_B/2$. 
Thus in the neighborhood of $H_B$, $\max_{k \in A_i\bs B_i} u_i(k,z_{-i}) < u_i(z)$ hence $R_{i,k}(z) < 0$ 
for all $k \in A_i\backslash B_i$. In particular, this holds if $z \in U_{\gamma_B}(H_B)$ (we omit the proof: easy but lengthy). 
It follows that if $z(t) \in U_{\gamma_B}(H_B)$, then by conditions (R3) and (Q2$'$), $a_i(t+1)\in B_i$ for each $i=1,2$.
\end{pf}

\bigskip
\begin{pf}[Proof of Claim \ref{Cm:1}]
Suppose that $\mathcal E_t$ holds and let $z(t)\in \mathcal U_{\gamma_B}(B)$. Then $a(t+1)\in B$  by Lemma \ref{Lem:U}.
We proceed by induction. Assume $a(t+1),\ldots,a(t+n)\in B$ for some $n\in \N^*$. Since $z(t)\in \mathcal U_{\gamma_B}(B)$,
\[
\sum_{a\in A\bs B} z_a(t+n) < \sum_{a\in A\bs B} z_a(t)<\gamma_B.
\]
Together with $\mathcal E_t$, this implies that $z(t+n)\in \mathcal U_{\gamma_B}(B)$. Consequently, by Lemma \ref{Lem:U}, $a(t+n+1)\in B$.
\end{pf}

\bigskip The proof of Claim 2 builds up on the proof of Theorem 2.1 of \citet{Hart01Gen}. It is different though, since we need to find the convergence rate of of the maximal regret {\it conditional} on a given initial history (in particular, on those where the past average play is close to a curb set), which \citet{Hart01Gen} do not provide.  So our result cannot be directly derived from their proof. 

For the proof of Claim \ref{Cm:2} we need the following lemmata.

\begin{lem}\label{Lem:K} 
Let $x_1,x_2,\ldots$ be a sequence of real random variables with $E[x_n|x_{n-1},\ldots,x_1]=0$ and $Var[x_n]\le \bar \sigma^2$ for all $n$. Then for every $\pi>0$ and every $m=1,2,\ldots$
\[
\Pr\left[\max_{n>m} \left|\frac{1}{n}\sum\nolimits_{k=m+1}^n x_k\right|\ge \frac{\bar \sigma}{\sqrt{m\pi}}\right]\le \pi.
\]
\end{lem}

\begin{pf}
H{\'a}jek-R{\'e}nyi inequality (e.g., \citet{Bullen98}) implies
\[
\Pr\left[\max_{m<k\le n} c_k\left|x_{m+1}+\ldots+x_k\right|\ge\eps\right]\le\frac{1}{\eps^2}\sum\nolimits_{k=m+1}^n c_k^2 Var[x_k].
\]
Using $c_k=1/k$ and $Var[x_k]\le \bar \sigma^2$, the right-hand side can be bounded as follows,
\[
\sum\nolimits_{k=m+1}^n c_k^2 Var[x_k]\le \bar \sigma^2 \sum\nolimits_{k=1}^{n-m} \frac{1}{(m+k)^2}\le \bar \sigma^2 \frac{1}{m}.
\]
Taking the limit $n\to\infty$ yields
\[
\Pr\left[\max_{k>m} \frac{1}{k}\left|x_{m+1}+\ldots+x_k\right|\ge\eps\right]\le\frac{\bar \sigma^2 }{m\eps^2}.
\]
The result is immediate by substitution $\pi=\bar \sigma^2/(m\eps^2)$.
\end{pf}

\bigskip
Define $\xi(1)=P_i(R_i(1))$ and for all $t=2,3,\ldots$ 
\begin{equation}\label{E:AQ3}
\xi(t)=tP_i(R_i(t))-(t-1)P_i(R_i(t-1)).
\end{equation}
\begin{lem}\label{Lem:L}
$\xi(t)$ is uniformly bounded and $E[\xi(t)| h(t-1)]\le C/t$ holds for some constant $C$ uniformly for all $t$.
\end{lem}

\begin{pf}
Let $x_0=R_i(t-1)$ and $x=R_i(t)$. Note that $R_i(t) = \frac{t-1}{t} R_i(t-1)+\frac{1}{t}r_i$, where $r_{i}=\left[u_i(k,a_{-i}(t))-u_i(a(t))\right]_{k\in A_i}$. Hence 
\begin{equation} \label{eq:xx0}
x-x_0=\frac{1}{t}(r_i-x_0).
\end{equation} 
The regret for an action is bounded by $2\bar U$ and the difference between two regret terms by $4 \bar U$. Thus, in sup norm, $|| r_i - x_0|| \leq 4 \bar U$ and $|| x-x_0 || \leq 4 \bar U / t$. 

Since $P_i$ is $C^2$, there exist constants $c$, $c'$ and $c''$ such that 
if $||y|| \leq   4 \bar U$,  
\[||P_i(y)|| \leq c, \quad ||\nabla P_i(y) \cdot  y|| \leq c' ||y||, \quad \mbox{ and }  ||y \cdot \nabla^2 P_i(y)  y|| \leq c'' ||y||^2.\] 
Moreover, $\xi(t) = P_i(x_0) + t(P_i(x)-P_i(x_0))$ hence $|\xi(t)| \leq c + t c' || x - x_0|| \leq c + 4\bar U c'$. Thus  $\xi(t)$ is uniformly bounded. 

We now show that $E[\xi(t)| h(t-1)]\le C/t$ for $C=8{\bar U}^2 c''$.  
By definition of $c''$ and Taylor-Lagrange theorem, 
\[P_i(x) \leq P_i(x_0) + \nabla P_i(x_0) \cdot (x-x_0) + \frac{1}{2} c'' ||x- x_0||^2.\]
  Using (\ref{eq:xx0}) we get: 
\[P_i(x) \leq \frac{t-1}{t}P_i(x_0)+\frac{1}{t}\left(P_i(x_0)-\nabla P_i(x_0) \cdot x_0\right)+\frac{1}{t}\nabla P_i(x_0) \cdot r_i(t) + \frac{C}{t^2}. \] 
Since $P_i$ is convex and $P_i(0)=0$, we have: 
\[P_i(x_0)  - \nabla P_i(x_0) \cdot x_0 = P_i(x_0) + \nabla P_i(x_0) \cdot (0-x_0)  \leq  P_i(0) = 0.\]
Therefore  
\[P_i(x) \leq \frac{t-1}{t}P_i(x_0)+ \frac{1}{t}\nabla P_i(x_0) \cdot r_i + \frac{C}{t^2}, \]
so that 
\[\xi(t) = tP_i(x) - (t-1)P_i(x_0) \leq  \nabla P_i(x_0) \cdot r_i + \frac{C}{t}. \]
To prove that $E[\xi(t)| h(t-1)]\le C/t$, it suffices to show that $E[ \nabla P_i(x_0) \cdot r_i | h(t-1)]=0$. To see this, note that 
$E[r_{i, k} | h(t-1)]= u_i(k, a_{-i}(t)) - u_i(q_i(t), a_{-i}(t))$ hence 
\[E[q_{i}(t) \cdot r_i | h(t-1)]= q_i(t) \cdot E[r_{i} | h(t-1)] = \sum_{k \in A_i} q_{i,k} [u_i(k, a_{-i}) - u_i(q_i, a_{-i})]= 0.\]
Since $q_i(t)$ is proportional to $\nabla P_i(x_0)$, the result follows. 
\end{pf}

\medskip
\begin{pf}[Proof of Claim \ref{Cm:2}] 
By \eqref{E:AQ3} we have
\begin{align*}
   P_i(R_i(t+n)) &=\frac{1}{t+n}\sum_{s=t+1}^{t+n} \xi(t)+\frac{t}{t+n}P_i(R_i(t)) \\
  &= \frac{1}{t+n}\sum_{s=t+1}^{t+n} \zeta(s)+ \frac{1}{t+n}\sum_{s=t+1}^{t+n} E[\xi(s)|h(s-1)]+\frac{t}{t+n}P_i(R_i(t)),
\end{align*}
where $\zeta(s)=\xi(s)-E[\xi(s)|h(s-1)]$. 
As we have assumed $z(t)\in\mathcal U_\gamma(B)$, we have
\[
\frac{t}{t+n}P_i(R_i(t))< \frac{t}{t+n} \gamma\le \gamma.
\]
Next, by Lemma \ref{Lem:L},
\[
\frac{1}{t+n}\sum_{s=t+1}^{t+n} E[\xi(s)|h(s-1)]\le \frac{1}{t+n}\sum_{s=t+1}^{t+n}\frac{C}{s}\le C\frac{\ln (t+n)-\ln t}{t+n}.
\]
Maximizing $\frac{\ln (t+n)-\ln t}{t+n}$ among all $n > 0$ yields $\frac{\ln (t+n)-\ln t}{t+n}\le (te)^{-1}$, hence
\[
\frac{1}{t+n}\sum_{s=t+1}^{t+n} E[\xi(s)|h(s-1)]\le \frac{C}{te}.
\]
Let $\bar\sigma^2$ be a bound on $Var[\zeta(s)]$ for all $s$ (this bound exists, since by Lemma \ref{Lem:L} variables $\xi(s)$ are uniformly bounded). Applying Lemma \ref{Lem:K}, we obtain that for every $\pi>0$ and every $t$,
\[
\Pr\left[\max_{n\in\N} \left|\frac{1}{t+n}\sum_{s=t+1}^{t+n} \zeta(s)\right|\ge \frac{\bar \sigma}{\sqrt{t\pi/2}}\right]\le \frac{\pi}{2}.
\]
Hence with probability at least $1-\pi/2$ the following holds for all $n\in\N^*$,
\[
   P_i(R_i(t+n))<\frac{\sqrt{2} \bar \sigma}{\sqrt{t \pi}}+ \frac{C}{te}+\gamma.
\]
and it holds for both $i=1,2$ simultaneously with probability at least $(1-\pi/2)^2\ge 1-\pi$.
Choosing $t_0=t_0(\pi,\gamma)$ so large that $\frac{\sqrt{2} \bar \sigma}{\sqrt{t_0 \pi}}+ \frac{C}{t_0e}\le \gamma_B-\gamma$, we obtain that event $\mathcal E_t$:
\[
   P_i(R_i(t+n))<\gamma_B \ \ \text{for each $i=1,2$ and all $n\in\N^*$}
\]
occurs with probability at least $(1-\pi/2)^2\ge 1-\pi$.
\end{pf}

\begin{spacing}{1}
\addtolength{\bibsep}{-3pt}

\end{spacing}


\begin{thebibliography}{47}
\providecommand{\natexlab}[1]{#1}
\providecommand{\url}[1]{\texttt{#1}}
\expandafter\ifx\csname urlstyle\endcsname\relax
  \providecommand{\doi}[1]{doi: #1}\else
  \providecommand{\doi}{doi: \begingroup \urlstyle{rm}\Url}\fi

\bibitem[Aubin and Celina(1984)]{Aubin84}
J.-P. Aubin and A.~Celina.
\newblock {\em Differential Inclusions}.
\newblock Springer, 1984.

\bibitem[Balkenborg et~al.(2012)Balkenborg, Hofbauer, and
  Kuzmics]{Balkenborg12}
D.~Balkenborg, J.~Hofbauer, and C.~Kuzmics.
\newblock Refined best reply correspondence and dynamics.
\newblock \emph{Theoretical Econ.}, forthcoming.

\bibitem[Basu and Weibull(1991)]{Basu91}
K.~Basu and J.~W. Weibull.
\newblock Strategy subsets closed under rational behavior.
\newblock \emph{Econ. Letters} 36 (1991), 141--146.

\bibitem[Bena{\"\i}m and Faure(2011)]{Benaim11}
M.~Bena{\"\i}m and M.~Faure.
\newblock Consistency of vanishing smooth fictitious play.
\newblock Working paper, available at http://arxiv.org/abs/1105.1690, 2011.

\bibitem[Bena{\"\i}m et~al.(2005)Bena{\"\i}m, Hofbauer, and Sorin]{Benaim05}
M.~Bena{\"\i}m, J.~Hofbauer, and S.~Sorin.
\newblock Stochastic approximations and differential inclusions.
\newblock \emph{SIAM J. Control and Optimization} 44 (2005), 328--348.

\bibitem[Bena{\"\i}m et~al.(2006)Bena{\"\i}m, Hofbauer, and Sorin]{Benaim06}
M.~Bena{\"\i}m, J.~Hofbauer, and S.~Sorin.
\newblock Stochastic approximations and differential inclusions. {Part II}:
  Applications.
\newblock \emph{Math. Operations Res.} 31 (2006), 673--695.

\bibitem[Berger(2005)]{Berger05}
U.~Berger.
\newblock Fictitious play in $2 \times n$ games.
\newblock \emph{J. Econ. Theory} 120 (2005), 139--154.

\bibitem[Berger(2007)]{Berger07}
U.~Berger.
\newblock Two more classes of games with the continuous-time fictitious play
  property.
\newblock \emph{Games Econ. Behav.} 60 (2007), 247--261.

\bibitem[Blackwell(1956)]{Blackwell56}
D.~Blackwell.
\newblock An analog of the minmax theorem for vector payoffs.
\newblock \emph{Pacific J. Math.} 6 (1956), 1--8.

\bibitem[Blum et~al.(2006)Blum, Even-Dar, and Ligett]{Blum06}
A.~Blum, E.~Even-Dar, and K.~Ligett.
\newblock Routing without regret: on convergence to Nash equilibria of
  regret-minimizing algorithms in routing games.
\newblock In \emph{Proceed. 25th Annual ACM Symposium on Principles
  of Distributed Computing}, pp.~45--52, 2006.

\bibitem[Brown(1951)]{Brown51}
G.~Brown.
\newblock Iterative solutions of games by fictitious play.
\newblock In T.~Koopmans (Ed.), \emph{Activity Analysis of Production and
  Allocation}, Vol.~13 of \emph{Cowles Commission Monograph}, pp.~374--376.
  New York: Wiley, 1951.
  
\bibitem[Bullen(1998)]{Bullen98}
P.~Bullen.
\newblock {\em A Dictionary of Inequalities}.
\newblock Addison Wesley Longman, Harlow, 1998.
  
\bibitem[Cesa-Bianchi and Lugosi(2003)]{CesaBianchi03}
N.~Cesa-Bianchi and G.~Lugosi.
\newblock Potential-based algorithms in on-line prediction and game theory.
\newblock \emph{Machine Learning} 51 (2003), 239--261.

\bibitem[Cesa-Bianchi and Lugosi(2006)]{CesaBianchi06}
N.~Cesa-Bianchi and G.~Lugosi.
\newblock \emph{Prediction, Learning, and Games}.
\newblock Cambridge Univ. Press, 2006.

\bibitem[Chen and Vaughan(2010)]{Chen10}
Y.~Chen and J.~W. Vaughan.
\newblock A new understanding of prediction markets via no-regret learning.
\newblock In \emph{Proceed. 11th ACM Conference on Electronic
  Commerce}, pp.~189--198, 2010.

\bibitem[Clemen and Winkler(2007)]{Clemen07}
R.~T. Clemen and R.~L. Winkler.
\newblock Aggregating probability distributions.
\newblock In W.~Edwards, R.~Miles, and D.~{von Winterfeldt} (Eds.),
  \emph{Advances Dec. Anal.}, pp.~154--176. Cambridge Univ. Press, 2007.

\bibitem[DeMarzo et~al.(2006)DeMarzo, Kremer, and Mansour]{DeMarzo06}
P.~DeMarzo, I.~Kremer, and Y.~Mansour.
\newblock Online trading algorithms and robust option pricing.
\newblock In \emph{Proceed. 38th Annual ACM Symposium on Theory of
  Computing}, pp.~477--486, 2006.

\bibitem[Faure and Roth(2010)]{Faure10}
M.~Faure and G.~Roth.
\newblock Stochastic approximations of set-valued dynamical systems: convergence with positive probability to an attractor.
\newblock \emph{Math. Operations Res.} 35 (2010), 624--640.

\bibitem[Foster and Vohra(1993)]{Foster93}
D.~Foster and R.~Vohra.
\newblock A randomization rule for selecting forecasts.
\newblock \emph{Operations Res.} 41 (1993), 704--709.

\bibitem[Foster and Vohra(1999)]{Foster99}
D.~Foster and R.~Vohra.
\newblock Regret in the online decision problem.
\newblock \emph{Games Econ. Behav.} 29 (1999), 7--35.

\bibitem[Freund and Schapire(1999)]{Freund99}
Y.~Freund and R.~Schapire.
\newblock Adaptive game playing using multiplicative weights.
\newblock \emph{Games Econ. Behav.} 29 (1999), 79--103.

\bibitem[Fudenberg and Levine(1995)]{Fudenberg95}
D.~Fudenberg and D.~Levine.
\newblock Consistency and cautious fictitious play.
\newblock \emph{J. Econ. Dynam. Control} 19 (1995), 1065--1089.

\bibitem[Gaunersdorfer and Hofbauer(1995)]{Gaunersdorfer95}
A.~Gaunersdorfer and J.~Hofbauer.
\newblock Fictitious play, {Shapley} polygons, and the replicator equation.
\newblock \emph{Games Econ. Behav.} 11 (1995), 279--303.

\bibitem[Gilboa and Matsui(1991)]{Gilboa91}
I.~Gilboa and A.~Matsui.
\newblock Social stability and equilibrium.
\newblock \emph{Econometrica} 59 (1991), 859--867.

\bibitem[Hannan(1957)]{Hannan57}
J.~Hannan.
\newblock Approximation to {Bayes} risk in repeated play.
\newblock In M.~Dresher, A.~W. Tucker, and P.~Wolfe (Eds.),
  \emph{Contributions to the Theory of Games, Vol. {III}}, Ann. Math. Stud. 39, pp.~ 97--139. Princeton Univ. Press, 1957.

\bibitem[Hart and Mas-Colell(2000)]{Hart00}
S.~Hart and A.~Mas-Colell.
\newblock A simple adaptive procedure leading to correlated equilibrium.
\newblock \emph{Econometrica} 68 (2000), 1127--1150.

\bibitem[Hart and Mas-Colell(2001)]{Hart01Gen}
S.~Hart and A.~Mas-Colell.
\newblock A general class of adaptive procedures.
\newblock \emph{J. Econ. Theory} 98 (2001), 26--54.

\bibitem[Hart and Mas-Colell(2003)]{Hart03}
S.~Hart and A.~Mas-Colell.
\newblock Continuous-time regret-based dynamics.
\newblock \emph{Games Econ. Behav.} 45 (2003), 375--394.

\bibitem[Hofbauer(1995)]{Hofbauer95}
J.~Hofbauer.
\newblock Stability for the best response dynamics.
\newblock Mimeo, 1995.

\bibitem[Hofbauer and Sandholm(2009)]{Hofbauer09JET}
J.~Hofbauer and W.~H. Sandholm.
\newblock Stable games and their dynamics.
\newblock \emph{J. Econ. Theory} 144 (2009), 1665--1693.

\bibitem[Hofbauer and Sigmund(1998)]{Hofbauer98}
J.~Hofbauer and K.~Sigmund.
\newblock \emph{Evolutionary Games and Population Dynamics}.
\newblock Cambridge Univ. Press, 1998.

\bibitem[Hofbauer and Sorin(2006)]{Hofbauer06}
J.~Hofbauer and S.~Sorin.
\newblock Best response dynamics for continuous zero-sum games.
\newblock \emph{Discrete and Continuous Dynamical Systems, Series B}, 6 (2006) 215--224.

\bibitem[Hofbauer et~al.(2009)Hofbauer, Sorin, and Viossat]{Hofbauer09}
J.~Hofbauer, S.~Sorin, and Y.~Viossat.
\newblock Time average replicator and best reply dynamics.
\newblock \emph{Math. Operations Res.} 34 (2009), 263--269.

\bibitem[Hofbauer and Weibull(1996)]{Hofbauer96}
J.~Hofbauer and J.~W. Weibull.
\newblock Evolutionary selection against dominated strategies.
\newblock {\em J. Econ. Theory} 71 (1996), 558--573.

\bibitem[Irani and Karlin(1996)]{Irani96}
S.~Irani and A.~Karlin.
\newblock On-line computation.
\newblock In D.~Hochbaum (Ed.), \emph{Approximation Algorithms for {NP}-Hard
  Problems}, pp.~521--564. Boston: PWS-Kent, 1996.

\bibitem[Krishna and Sj{\"o}str{\"o}m(1998)]{Krishna98}
V.~Krishna and T.~Sj{\"o}str{\"o}m.
\newblock On the convergence of fictitious play.
\newblock \emph{Math. Operations Res.} 23 (1998), 479--511.

\bibitem[Larrick and Soll(2006)]{Larrick06}
R.~P. Larrick and J.~B. Soll.
\newblock Intuitions about combining opinions: misappreciation of the averaging
  principle.
\newblock \emph{Manage. Sci.} 52 (2006), 111--127.

\bibitem[Lehrer(2003)]{Lehrer03}
E.~Lehrer.
\newblock A wide range no-regret theorem.
\newblock \emph{Games Econ. Behav.} 42 (2003), 101--115.

\bibitem[Littlestone and Warmuth(1994)]{Littlestone94}
N.~Littlestone and M.~Warmuth.
\newblock The weighted majority algorithm.
\newblock \emph{Information and Computation} 108 (1994), 212--261.

\bibitem[Mansour(2010)]{Mansour10}
Y.~Mansour.
\newblock Regret minimization and job scheduling.
\newblock In \emph{Proceed. 36th Conference on Current Trends in
  Theory and Practice of Computer Science}, pp.~71--76. Springer, 2010.

\bibitem[Matsui(1992)]{Matsui92}
A.~Matsui.
\newblock Best response dynamics and socially stable strategies.
\newblock \emph{J. Econ. Theory} 57 (1992), 343--362.

\bibitem[Monderer et~al.(1997)Monderer, Samet, and Sela]{Monderer97}
D.~Monderer, D.~Samet, and A.~Sela.
\newblock Belief affirming in learning processes.
\newblock \emph{J. Econ. Theory} 73 (1997), 438--452.

\bibitem[Moulin and Vial(1978)]{Moulin78}
H.~Moulin and J.~P. Vial.
\newblock Strategically zero-sum games: the class of games whose completely
  mixed equilibria cannot be improved upon.
\newblock \emph{Int. J. Game Theory} 7 (1978), 201--221.

\bibitem[Selten(1995)]{Selten95}
R.~Selten.
\newblock An axiomatic theory of a risk dominance measure for bipolar games
  with linear incentives.
\newblock \emph{Games Econ. Behav.} 8 (1995), 213--263.

\bibitem[Shapley(1964)]{Shapley64}
L.~S. Shapley.
\newblock Some topics in two-person games.
\newblock In M.~Dresher, L.~S. Shapley, and A.~W. Tucker, editors,
  \emph{Advances in Game Theory}, pp.~1--28. Princeton Univ. Press,
  1964.

\bibitem[Sparrow et~al.(2008)Sparrow, van Strien, and Harris]{Sparrow08}
C.~Sparrow, S.~van Strien, and C.~Harris.
\newblock Fictitious play in $3\times 3$ games: The transition between periodic
  and chaotic behaviour.
\newblock \emph{Games Econ. Behav.} 63 (2008), 259--291.

\bibitem[Timmerman(2006)]{Timmerman06}
A.~Timmerman.
\newblock Forecast combinations.
\newblock In G.~Elliott, C.~W. Granger, and A.~Timmermann (Eds.),
  \emph{Handbook of Economic Forecasting}. Elsevier, 2006.

\bibitem[Tirole(1988)]{Tirole88}
J.~Tirole.
\newblock \emph{The Theory of Industrial Organization}.
\newblock MIT Press, 1988.

\bibitem[Weibull(1995)]{Weibull95}
J.~W.~Weibull.
\newblock {\em Evolutionary Game Theory}.
\newblock Cambridge Univ. Press, 1995.


\bibitem[Yanovskaya(1968)]{Yanovskaya68}
E.~Yanovskaya.
\newblock Equilibrium points in polymatrix games (in {Russian}).
\newblock \emph{Litovskii Matematicheskii Sbornik} 8 (1968), 381--384.

\bibitem[Young(1993)]{Young93E}
H.~P. Young.
\newblock The evolution of conventions.
\newblock \emph{Econometrica} 61 (1993), 57--84.

\bibitem[Young(2004)]{Young04}
H.~P. Young.
\newblock \emph{Strategic Learning and Its Limits}.
\newblock Oxford Univ. Press, 2004.

\end{thebibliography}
\end{document}